\documentclass[journal]{new-aiaa}
\usepackage[utf8]{inputenc}
\usepackage{textcomp}
\usepackage{booktabs}
\usepackage{multirow}

\usepackage{graphicx}
\usepackage{subcaption}
\usepackage{amsmath}
\usepackage[version=4]{mhchem}
\usepackage{siunitx}
\usepackage{longtable,tabularx}
\usepackage{xurl}
\usepackage{xcolor}
\newcolumntype{C}[1]{>{\PreserveBackslash\centering}p{#1}}
\newcolumntype{R}[1]{>{\PreserveBackslash\raggedleft}p{#1}}
\newcolumntype{L}[1]{>{\PreserveBackslash\raggedright}p{#1}}

\title{\color{black}System Architecting for GEO Communication Satellite Considering On-Orbit Refueling}

\author{Jaewoo Kim \footnote{PhD Candidate, Department of Aerospace Engineering, 291 Daehak-Ro.}}
\affil{Korea Advanced Institute of Science and Technology, Daejeon 34141, Republic of Korea}
\author{Jaemyung Ahn \footnote{Professor of Aerospace Engineering, 291 Daehak Ro; jaemyung.ahn@kaist.ac.kr. AIAA Associate Fellow (Corresponding Author).}}
\affil{Korea Advanced Institute of Science and Technology, Daejeon 34141, Republic of Korea}

\begin{document}

\maketitle

\begin{abstract}
    This paper introduces the problem of selecting a satellite system architecture considering commercial on-orbit refueling (OOR). {\color{black}This problem answers two questions: ``What design lifetime should the satellite have?'' and ``How much propellant should be carried at launch?''} We formulate {\color{black}this as a mathematical optimization problem by adopting design lifetime and initial propellant mass} as design variables and considering two objective functions to balance the returns and risks. {\color{black}To solve this problem, we develop a surrogate model-based framework grounded in a satellite lifecycle simulation. The framework captures} various uncertainties and operational flexibility, and integrates a modified satellite sizing and cost model by adjusting traditional models with OOR. {\color{black} Based on the developed framework, we conduct a case study of GEO communication satellites to examine current target service performance and explore the potential of a new system architecture that diverges from traditional design trends.}
\end{abstract}

\section*{Nomenclature}
{\renewcommand\arraystretch{1.0}
\noindent\begin{longtable*}{@{}l @{\quad=\quad} l@{}}
$a_{\text{prop}},b_{\text{prop}}$ & coefficients for propulsion subsystem mass, kg$^{1/3}$, kg\\
$C_{\text{IOC}}$ & initial operational capability cost, \$M\\
$C_{\text{lau}}$ & launch cost, \$M\\
$C_{\text{oor}}$ & On-Orbit Refueling (OOR) service cost, \$M\\
$C_{\text{op}}$ & satellite operational cost, \$M\\
$C_{\text{sat}}$ & satellite acquisition cost, \$M\\
\color{black}$C_{\text{serv}}$ & \color{black}cost of integrating servicing instruments, \$M\\
$c_{\text{lau}}$ & specific launch cost, \$M/kg\\
\color{black}$c_{\text{oor},f}$ & fixed cost of OOR service, \$M\\
\color{black}$c_{\text{oor},v}$ & variable cost coefficient of OOR service, \$M/kg\\
$c_{\text{sat}}$ & specific satellite acquisition cost, \$M/kg\\
$I_{\text{sp}}$ & specific impulse of the propulsion subsystem, sec\\
$M_{\text{oor}}$ & OOR service capacity, kg\\
$m_{\text{adcs}}$ & attitude determination and control subsystem mass of a satellite, kg\\
$m_{\text{base}}$ & base mass of a satellite, kg\\
$m_{\text{dry}}$ & dry mass of a satellite, kg\\
$m_{\text{ps}}$ & propulsion subsystem mass of a satellite, kg\\
$m_{\text{ref}}$ & reference mass, kg\\
$m_{\text{serv}}$ & servicing interface mass of a satellite, kg\\
\color{black}$m_{\text{serv}}$ & \color{black}mass of servicing instruments, kg\\
$m_{\text{str}}$ & structure mass of a satellite, kg\\
$m_{\text{wet}}$ & wet mass of a satellite, kg\\
$m_{p,\text{des}}$ & design propellant mass of a satellite, kg\\
$m_{p,\text{ioa}}$ & propellant mass required for initial orbit acquisition, kg\\
$m_{p,\text{oor}}$ & propellant amount to be replenished, kg\\
\color{black}$p_{\text{lau}}$ & launch failure rate\\
\color{black}$p_{\text{oor}}$ & OOR service failure rate\\
$R$ & satellite operational revenue, \$M/yr\\
$R_0$ & initial operational revenue, \$M/yr\\
$Rel$ & reliability function of a satellite\\
$Rel_{\text{ref}}$ & reliability function of a satellite with the reference design lifetime\\
$r_f$ & \color{black} continuously compounded risk free rate\\
$T$ & time in continuous frame, yr\\
$T_{\text{life}}$, $t_{\text{life}}$ & design lifetime of satellite in year and time step, yr, -\\
\color{black}$T_{\text{oor}}$, $t_{\text{oor}}$ & time to OOR service in year and time step, yr, -\\
$T_{\text{ref}}$ & reference design lifetime, yr\\
\color{black}$T_{\text{rep}}$, $t_{\text{rep}}$ & time to replacement in year and time step, yr, -\\
$T_{\text{sim}}$, $t_{\text{sim}}$ & simulation time in year and time step, yr, -\\
$t$ & time in discretized frame ($=T/\Delta T$)\\
$\Delta T$ & time step size in year, yr\\
$\Delta V_{\text{err}}$ & $\Delta V$ required for correcting the injection error of the launch vehicle, m/s\\
$\Delta V_{\text{ioa}}$ & $\Delta V$ required for initial orbit acquisition, m/s\\
$\Delta V_{\text{ot}}$, $\Delta V_{\text{ot},\text{ideal}}$ & real and ideal $\Delta V$ required for transferring from the initial orbit to the target orbit, m/s, m/s\\
$\Delta V_{\text{stk},\text{yr}}$ & {\color{black}a}nnual $\Delta V$ requirement for station-keeping, m/s\\
\color{black}$\Delta V_{\text{stk}}$ & $\Delta V$ requirement for station-keeping per time step, m/s\\
\color{black}$\alpha_{\text{adcs}}$ & attitude determination and control subsystem mass ratio\\
\color{black}$\alpha_{\text{ins}}$ & insurance cost ratio\\
\color{black}$\alpha_{\text{op}}$ & operational cost ratio\\
\color{black}$\alpha_{\text{str}}$ & structural mass ratio\\
\color{black}$\alpha_{\text{rel}},\beta_{\text{rel},1},\beta_{\text{rel},2},\theta_{\text{rel},1},\theta_{\text{rel},2}$ & parameters for reliability function, -, -, -, yr, yr\\
\color{black}$\theta_{\text{obs}}$ & parameter for technology obsolescence, yr\\
$\kappa$ & mass growth rate of a satellite to the design lifetime\\
$\mu_{\text{mar}}$ & market drift rate, yr$^{-1}$\\
$\sigma_{\text{mar}}$ & market volatility, yr$^{-1/2}$\\
$\sigma_{\text{oi}}$ & parameter of $\Delta V$ distribution for correcting orbit injection error\\
$\Phi$ & cumulative distribution function of the standard normal distribution\\
$\phi_{\text{mar}}$ & market volatility factor in operational revenue\\
$\phi_{\text{obs}}$ & technology obsolescence factor in operational revenue\\
\color{black}$\mathcal{A}$ & \color{black}experiment design variables set\\
$\mathcal{H}$ & half-normal distribution with zero as the lower bound\\
\color{black}$\mathcal{X}$ & \color{black}feasible set of design variables\\
\color{black}DCF &\color{black}discounted cash flow\\
\color{black}FCF &\color{black}free cash flow\\
\color{black}MOO &\color{black}multi-objective optimization\\
\color{black}NPV &\color{black}net present value \\
\color{black}OOR &\color{black}on-orbit refueling \\
\color{black}OOS &\color{black}on-orbit servicing \\
\end{longtable*}}

\section{Introduction}\label{sec: section 1}
\lettrine{S}{PACE} systems remain among the most expensive and complex, even in this era of abundant advanced technologies. However, only a limited number of these systems, such as the International Space Station and the Hubble Space Telescope, have benefited from post-launch servicing (e.g., maintenance, repair, or replenishment)---a routine practice for complex ground systems. {\color{black}The primary reasons for this limitation are high costs and risks associated with human-involved on-orbit servicing (OOS) and the expenses of launches.} Consequently, most space systems are equipped with all the necessary resources for their entire operation at launch. This approach, while necessary, has exacerbated the impact of failures, thereby placing a premium on system reliability. The increased focus on reliability and the significant price of launch vehicles {\color{black}has} led to escalating expenses and complexity in satellite systems. This trend, characteristic of the risk-averse space industry, has contributed to the gradual extension of the design lifetimes of space systems.

{\color{black}However, the space industry is undergoing a significant transformation with the emergence of new commercial entities. These entities are driven by advancements in launch services, which have significantly reduced launch costs, fostering innovative concepts that enhance the value of space systems.} A notable innovation is the on-orbit refueling (OOR), a subset of OOS \cite{li2019orbit, arney2021orbit, davis2019orbit, cavaciuti2022space}. This innovation, supported by advancements in robotics, autonomous control, and servicing architecture design, primarily serves as a life-extension solution for existing assets of substantial value, such as the GEO communication satellites. This approach does not require significant modifications to existing design and operational paradigms. Instead, it introduces supplementary procedures after {\color{black}the initial propellant load is depleted}, promising profitability through the extended {\color{black}utilization} of these satellites. Consequently, the commercialization of OOR appears imminent as stakeholders in the space industry seek to capitalize on its potential. For example, Orbit Fab plans to provide in-orbit hydrazine for GEO satellites weighing up to 100 kilograms at a cost of {\color{black}\$20 million by 2025 \cite{orbitfab}.}

On the other hand, OOR also introduces significant changes to satellite system design---the separation of the design lifetime and the design propellant mass of a satellite. Traditionally, {\color{black}a satellite carries all the propellant needed throughout its mission duration.} However, with the OOR, satellites can carry a smaller initial propellant amount and rely on the refueling service in the future. This raises a question about the traditional approach: ``Will it remain the best option to load all the propellant required for the lifetime at launch?{\color{black} If not, how much should be loaded initially, and when---and in what quantity---should it be replenished?''} These questions lead to studying the potential of architectural changes, which is essential not only for the profitability {\color{black}of market players} but also for the long-term sustainability of the overall space industry \cite{hastings2016when}. Answering these questions requires consideration of {\color{black}satellites' behavior, interactions with the operational environment, and how these factors influence their overall value. A comprehensive evaluation of the satellite's entire lifecycle is essential to analyze this dynamic and reach the best decision.}

The design lifetime of a system---the period during which the system operates as intended---appears as a requirement and derives other lower-level requirements, such as the reliability of each subsystem and the amount of each consumable. The satellite's design lifetime has been a research topic{\color{black}, a }crucial element of satellite architecture. Saleh et al. \cite{saleh2002spacecraft} analyzed the impact of the satellite design lifetime on the sizing and cost of satellites and provided models for these as functions of the design lifetime{\color{black}.}

In a commercial context, economic considerations profoundly influence the decision-making process regarding spacecraft design lifetime, as explored in a series of papers {\color{black}\cite{saleh2004weaving, saleh2006reduce}. These studies examined how design lifetime impacts the value of satellite systems, considering environmental factors (e.g., market and technological advancements) from multiple perspectives, including those of satellite operators, manufacturers, society, and customers. Overall, these studies imply that as situations become increasingly competitive and uncertain, a simplistic, cost-driven approach leads to misguided solutions that expose stakeholders to greater risks.}

In addition to these advancements in analyzing the economic implications of satellite design lifetime, the value analysis of OOS in space systems has been extensively studied since the 1970s. The research group led by Dr. Hastings is a significant contributor to this field. They focused on the strategic advantage of flexibility provided by OOS, employing quantitative approaches to analyze its value in various OOS scenarios. Their extensive work, especially in the 2000s, is well-documented in various journal articles \cite{saleh2003flexibility, saleh2002space, lamassoure2002space, joppin2006orbit, long2007orbit} and theses \cite{lamassoure2001framework, saleh2002weaving, joppin2004orbit, long2005framework, richards2006orbit}, offering a comprehensive view of the strategic implications and benefits of OOS. A key contribution of these works is the incorporation of the operational flexibility provided by OOS into the evaluation framework through various techniques (e.g., decision tree analysis and real options \cite{mun2016real}). Using these techniques, the authors investigated the value offered by different OOS types---such as life-extension and upgrades---and examined how OOS could lead to changes in satellite system architecture, specifically in the design lifetime of a satellite.

Recognizing the flexibility OOS provides has prompted extensive research into its value across various operational scenarios, employing diverse methodologies. For example, Yao et al. \cite{yao2013orbit} approached OOS from a design optimization perspective, considering both service providers and clients. Their methodology factored in design lifetime, life-extension length, number of servicing operations before refueling, and the orbit radii of parking and depots. They utilized probability and evidence theories to address uncertainties in the operational environment. The study employed dynamic programming to select the optimal servicing activity based on modeled uncertainties. The authors decomposed the optimization problem into sub-problems to manage the computational complexity, enabling efficient analysis of a hypothetical GEO communication satellite case.

{\color{black}More recently, Liu et al. \cite{liu2021economic} focused on the economic aspects of OOS. They developed a high-fidelity cost-and-benefit model, incorporating insurance, taxes, depreciation, and other realistic economic parameters. With this model, they examined the MEV-1 mission \cite{henry2020intelsat}, a life-extension mission for Intelsat 901.} Their analysis revealed that the MEV-1's service pricing was at the lower end of the feasible range, likely to attract customers in a nascent market. They also suggested that the distribution of value and risk associated with OOS will be crucial for the future space industry.

Building on these academic efforts, we introduce a new problem in satellite system architecting. In this problem, we decouple the design lifetime from the propellant mass to account for an expanded range of satellite system architectural solutions that OOR commercialization drives---a perspective not covered in existing works. To address this problem while accommodating operational flexibility and various uncertainties, we propose a satellite lifecycle simulation for evaluating each architectural solution and introduce a solution procedure assisted by surrogate models.

{\color{black}
Throughout this research, we focus on GEO communication satellites for three reasons. First, GEO satellites represent one of the most successful commercial space systems to date. Second, their design has continually improved capacity and longevity, yielding reliable platforms exceeding six metric tons in wet mass \cite{krebs}, consistent with the aforementioned trends in space systems design. Third, the substantial costs and potential revenues over a GEO satellite’s lifecycle have made them primary candidates for OOR by emerging market players---e.g., Orbit Fab \cite{orbitfab}---unlike smaller satellites in low Earth orbit.

The rest of this paper is organized into five sections. Section \ref{sec: section 2} introduces the decision-making problem of interest and provides its mathematical representation. Section \ref{sec: section 3} constructs the satellite lifecycle simulation, which serves as the baseline for the solution procedure. Section \ref{sec: section 4} introduces the solution procedure, assisted by the constructed simulation and surrogate models. Section \ref{sec: section 5} presents two case studies---one for a GEO communication satellite with chemical propulsion and the other for a GEO communication satellite with electric propulsion---to validate the developed framework and offer meaningful insights. Finally, Section \ref{sec: section 6} concludes the paper, summarizing its contents and contributions and proposing future works.}

\section{Problem Description}\label{sec: section 2}
In a traditional approach, the design lifetime of a satellite, which is constrained by the planned operational time and technological factors, primarily influences its design propellant mass. Consider, for instance, a GEO communication satellite intended for a 15-year operational period. In such cases, {\color{black}this} satellite is launched with enough propellant to support all the projected maneuvers during its mission. This covers initial orbit insertion errors and station-keeping maneuvers, which typically require a $\Delta V$ of around 50 m/s annually. Within this design scheme, the designer's choice of propellant mass is highly constrained. Without {\color{black}an} additional supply strategy except the initial loading, the degree of conservatism becomes the only variable under the designer's control. This degree of conservatism typically appears in the design as a margin, providing a buffer to ensure mission success under potential deviations in mission environments.

However, while this traditional approach is time-tested, it limits the design space of satellites in the context of OOR. OOR introduces a paradigm shift by allowing propellant replenishment after launch, expanding the design space beyond those dictated solely by the loaded propellant at launch. In this new context, the propellant mass is no longer solely tied to the spacecraft's design lifetime. Instead, it becomes an independent variable determined by a comprehensive assessment of how it impacts the overall system value. This change enables a more adaptable satellite system architecture, potentially enhancing the system's value. A careful approach is essential in the initial design phase for the designer to realize this potential.

This approach involves {\color{black}a} trade-off between efficiency improvement and {\color{black}increased complexity.} On the one hand, reducing the initial propellant mass leads to a smaller tank than traditional designs, {\color{black}thereby} decreasing the satellite's overall structural mass. {\color{black}This reduction in structural mass may relax the performance requirements on the actuators, such as reaction wheels. Overall, reducing the initial propellant mass can lower mission costs, including development, manufacturing, and launch.} Additionally, the reduced overall expenditure on system acquisition mitigates {\color{black}the} loss in mission failures.

On the other hand, integrating OOR into the satellite's operational scenario introduces new complexities. This integration not only incurs additional servicing costs and introduces a potential new source of failure but also requires accommodating the service interface for OOR, which increases the satellite's size. The cumulative effect of these factors could offset the benefits derived from reduced propellant mass. Consequently, while OOR offers opportunities for cost reduction and increased flexibility, it also presents challenges that necessitate cautious management during the initial design phase.

{\color{black}A well-structured methodology grounded in quantitative analysis is required to evaluate these trade-offs and identify an efficient solution systematically.} For this purpose, {\color{black}we first mathematically encode} this decision-making context into an optimization problem. {\color{black}We introduce the satellite’s design lifetime and propellant mass as independent design variables.} Each choice of design variables reflects the designer's intention about ``how long the system should operate'' and ``how much propellant {\color{black}to load at launch,'' with each pairing representing a distinct design solution.}

{\color{black}The optimal values of the design variables depend on the chosen criteria. One possible approach is to use the monetary value generated by the satellite as the criterion,} which fits a commercial context. {\color{black}Additionally, in many cases, other forms of value can be converted into monetary terms using an appropriate model. With this approach}, the value unfolds as {\color{black}a} satellite operates throughout its lifetime, fulfilling its intended functions. This inherently uncertain process leads satellite operators {\color{black}to maximize profits while minimizing risks. We adopt} the metrics used for Markowitz’s mean-variance method in portfolio selection theory \cite{markowitz1952portfolio, sharpe1966mutual} with necessary modifications.

{\color{black}We formulate} the optimization problem for determining the optimal satellite system architecture {\color{black}as follows:}
\newline\newline\noindent
($\mathbf{P}_{\text{O}}$) Optimization of Satellite System Architecture
\begin{equation}
    \underset{\mathbf{x}}{\text{maximize}}\ \mathbf{J}=\left[E_V,\ E_V/\sigma_V\right]{\color{black},}
\end{equation}
subject to
\begin{equation}
    \mathbf{x}=\left[T_{\text{life}},\ m_{p,\text{des}}\right]\in {\color{black}\mathcal{X}}.
\end{equation}
Here, $T_{\text{life}}$, $m_{p,\text{des}}$, {\color{black}and $\mathcal{X}$ are} design lifetime, propellant mass, {\color{black}the feasible set of design variables,} respectively, and $V$ denotes the {\color{black}net present value (NPV)} representing the present value of cash flows throughout the project. {\color{black}The first objective, $E_V$, is the expected NPV, reflecting profit maximization. The second objective, the ratio of expected NPV versus its standard deviation, $E_V/\sigma_V$, is analogous to the Sharpe ratio \cite{sharpe1966mutual}, representing the profit per unit risk.}

This multi-objective optimization (MOO) problem has a family of efficient solutions, {\color{black}known as the Pareto front}. When solved {\color{black}by} an appropriate MOO algorithm, the problem yields a set of non-dominated solutions. {\color{black}Here, a non-dominated solution is one that is not inferior to any other solution with respect to all objective functions} \cite{marler2004survey}. In this formulation, each non-dominated solution offers a unique expected return and a specific expected return per risk, reflecting an efficient architectural solution.

\section{Satellite Lifecycle Simulation}\label{sec: section 3}
Traditionally, the NPV of an engineering project is evaluated using the discounted cash flow (DCF) method. This method projects future cash flows and converts them into present values by discounting with a {\color{black}discount} rate \cite{newnan2017engineering}. NPV is the sum of these present values ({\color{black}i.e.,} discounted cash flows). The chosen discount rate reflects the time value of money and {\color{black}risk} considerations associated with the projected cash flows. The discount rate is typically determined using established models (e.g., {\color{black}weighted average cost of capital} (WACC) \cite{miles1980weighted}) supported by empirical data. The DCF method is generally used in various engineering fields due to its simplicity. However, the method encounters significant challenges in newly developed concepts like space systems with OOR. The difficulty stems from the novelty and complexity of these systems, which {\color{black}complicate} the determination of an appropriate discount rate. Additionally, the DCF method often struggles to incorporate operational flexibility and dynamic risk structures.

To overcome these limitations, we propose a simulation-based approach. This simulation integrates key factors influencing satellite profitability, effectively representing operational scenarios. The methodology facilitates NPV estimations under various scenarios, each influenced by selected decision variables and parameters. This approach allows for a more realistic financial performance assessment reflecting the satellites' unique characteristics. By simulating the satellite's entire lifecycle, the simulation yields valuable insights into how different design lifetimes and propellant masses impact the financial performance of the satellite. The rest of this section details the simulation and its components.

\subsection{Satellite Sizing and Cost Elements}
Saleh et al. \cite{saleh2002spacecraft} established a mass estimating relationship (MER) that scales a satellite's mass as a function of its design lifetime. This MER has been used in various OOS evaluation frameworks until recently \cite{yao2013orbit, liu2021economic}. This approach assumes a uniform design lifetime across all subsystems and derives mass and cost elements based on this assumption.

{\color{black}We modify} this model by incorporating the propellant mass as an independent variable, yielding a new MER that integrates design lifetime and propellant mass as independent parameters. Subsequently, cost elements are derived from this revised MER.

The wet mass ($m_{\text{wet}}$) is the sum of the dry mass ($m_{\text{dry}}$) and the design propellant mass ($m_{p,\text{des}}$):
\begin{equation}
    m_{\text{wet}} = m_{\text{dry}} + m_{p,\text{des}}.
    \label{eqn: MER1}
\end{equation}
Here, the dry mass is decomposed into subsystem masses. We consider subsystems that are significantly influenced in size by the design propellant mass---the propulsion subsystem ($m_{\text{ps}}$), the structural subsystem ($m_{\text{str}}$), and the attitude determination and control subsystem ($m_{\text{adcs}}$)---as individual elements. The mass increase due to the servicing interface ($m_{\text{serv}}$) is considered a distinct element. {\color{black}All other subsystems are combined into a single integrated base mass ($m_{\mathrm{base}}$), while the payload mass ($m_{\mathrm{pl}}$) is treated separately.
\begin{equation}
    m_{\text{dry}} = m_{\text{base}} + m_{\text{pl}}+ m_{\text{ps}} + m_{\text{str}} + m_{\text{adcs}} + m_{\text{serv}}.
    \label{eqn: MER2}
\end{equation}
Each element of the dry mass is modeled as follows:
\begin{equation}
    m_{y} = m_{y,\text{ref}} \frac{1 + \kappa \left(T_{\text{life}} - 3\right)}{1 + \kappa \left(T_{\text{ref}} - 3\right)},\quad \text{for } y=\text{base},\ \text{pl},
    \label{eqn: MER3}
\end{equation}
\begin{equation}
    m_{\text{ps}} = a_{\text{prop}}  m_{p,\text{des}}^{2/3} + b_{\text{prop}},
    \label{eqn: MER4}
\end{equation}
\begin{equation}\label{eqn: MER5}
    m_{y} = \alpha_{y} m_{\text{dry}},\quad \text{for } y=\text{adcs},\ \text{str}.
\end{equation}
Here, Eqs. \eqref{eqn: MER3} and \eqref{eqn: MER4} adopt the relations referred to \cite{saleh2002spacecraft}. In particular, Eq. \eqref{eqn: MER3} implies $\kappa$\% increase in mass for each additional year of design lifetime beyond three years. We then model $m_{\text{adcs}}$ and $m_{\text{str}}$ as proportional to the dry mass, as shown in Eq. \eqref{eqn: MER5}.}

We introduce cost models that minimize dependency on additional parameters and reduce modeling complexity. Three main events incur costs: {\color{black}satellite replacement, satellite operations,} and OOR. Satellite replacement incurs the cost required to achieve the initial operating capability ($C_{\text{IOC}}$){\color{black}, which is defined as follows:
\begin{equation}\label{eqn: CER initial operational capacity}
    C_{\text{IOC}} = \left(1+\alpha_{\text{ins}}\right)  C_{\text{sat}} + C_{\text{lau}}.
\end{equation}
Here, $\alpha_{\text{ins}}$ is the insurance cost ratio, $C_{\text{sat}}$ the satellite acquisition cost, and $C_{\text{lau}}$ is the launch cost. The acquisition cost $C_{\text{sat}}$ can be estimated using various available parametric models; we adopt the unmanned space vehicle cost model (USCM) 8 \cite{wertz2011space}:
\begin{equation}\label{eqn: CER satellite}
    C_{\text{sat}} = 1.124\times 1.234\times \left(283.5\left(m_{\text{dry}}-m_{\text{pl}}-m_{\text{serv}}\right)^{0.716}+189m_{\text{pl}}+C_{\text{serv}}\right)\times\frac{\text{CPI}_{2025}}{\text{CPI}_{2010}}\times\frac{1}{1000},
\end{equation}
where $C_{\mathrm{serv}}$ is the cost of integrating servicing instruments to enable OOR (e.g. servicing interface), $\text{CPI}_{2010}$ and $\text{CPI}_{\text{2025}}$ are consumer price index (CPI) of 2010 and 2025 in the US \cite{cpi}, respectively. Also, we assume launch cost as proportional to the wet mass:
\begin{equation}\label{eqn: CER launch}
    C_{\text{lau}} = c_{\text{lau}}  m_{\text{wet}}.
\end{equation}}

We model the service price of OOR as a function of the propellant mass to be delivered to the client satellite, adopting a servicing infrastructure-neutral cost model. This approach does not rely on specific assumptions about the servicing infrastructure. Instead, it reflects the pricing policy of the service provider{\color{black}. Under a linear pricing policy, the service cost is expressed as follows:
\begin{equation}
    C_{\text{oor}} = c_{\text{oor},v}  m_{p,\text{oor}} + c_{\text{oor},f}.
    \label{eqn: CER3}
\end{equation}}

Lastly, the satellite operation cost ($C_{\text{op}}$) is modeled to be proportional to the revenue provided by the satellite operation ($R$):{\color{black}
\begin{equation}
    C_{\text{op}} = \alpha_{\text{op}} R.
    \label{eqn: CER4}
\end{equation}}

\subsection{Construction of the Satellite Lifecycle Simulation}
We {\color{black}construct} a computer simulation describing the states and actions of a satellite throughout its lifecycle. A given time horizon ($T_\text{sim}$) for the simulation is divided into a constant interval ($\Delta T$), yielding discretized time steps ($t = 0,\dots, t_{\text{sim}}; t_{\text{sim}} = T_{\text{sim}}/\Delta T$).
\begin{figure}[hbt!]
	\centering
	\includegraphics[page=1,width=1.0\columnwidth]{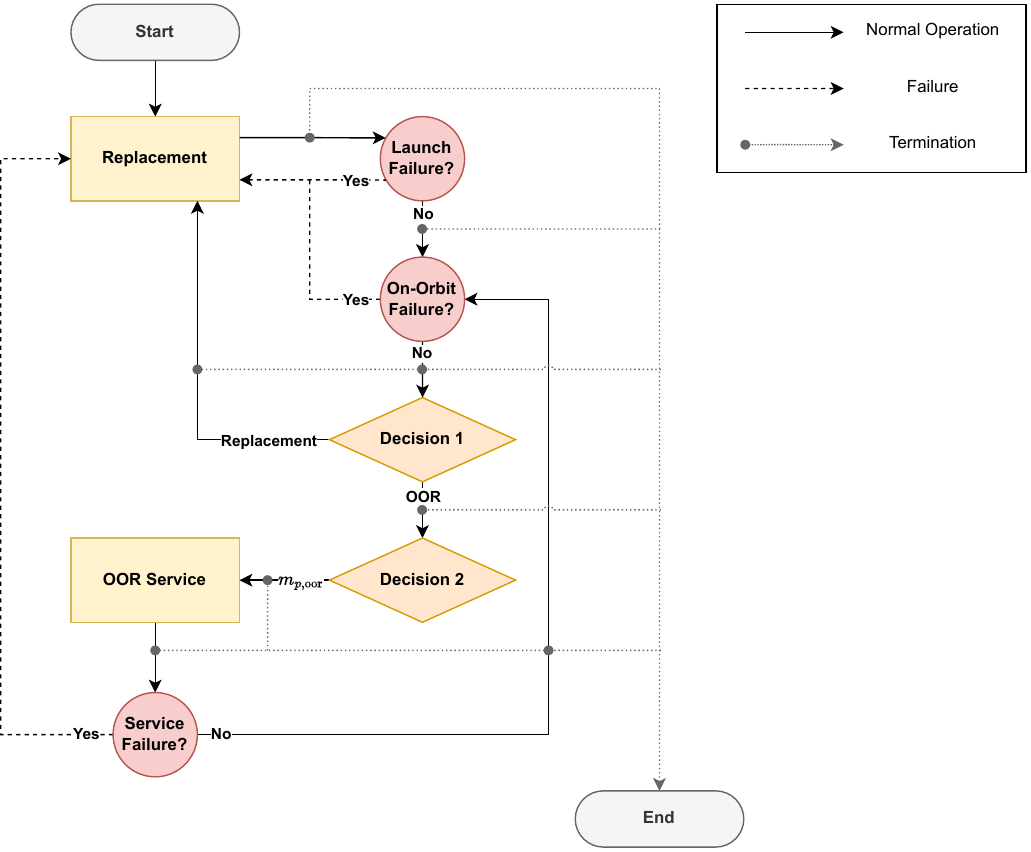}
	\caption{Flow diagram of the satellite lifecycle simulation{\color{black}.}}
    \label{fig: flow diagram}
\end{figure}

The assumptions made for the simulation are as follows:
\begin{enumerate}
\item Cost Incurrence: Cost is incurred at the time of the corresponding event. For example, a new satellite launch incurs the replacement cost, implying we focus on the launch event rather than the satellite manufacturing and launch contracting process.
\item Revenue Generation: Revenue occurs at time step $\left(t+1\right)$ if {\color{black}a} satellite is operational {\color{black}from} time step $t$ {\color{black}to} time step $\left(t+1\right)$ without failure.
\item Independence of Uncertainties: Different types of uncertainties, such as market {\color{black}fluctuation} and technological advancement, are assumed to be independent.
\item Replacement Failure: This failure type includes catastrophic launch failure and {\color{black}insufficient} propellant for initial orbital transfer. No rescue strategy for satellites failing to reach the intended mission orbit is assumed.
\item Insurance Coverage: In the event of failures during replacement or OOR service, associated costs are assumed to be covered by insurance. 
\item Technological Update: {\color{black}Each newly launched satellite incorporates the most up-to-date technology available at launch.}
\item Exclusion of Project Abandonment: The option to abandon the project is not considered.
\end{enumerate}

The simulation progresses by transitioning from one event to another. The event includes the \textit{Start} and the \textit{End} of the simulation, \textit{Decision}, \textit{Action}, and \textit{Failure}. Figure \ref{fig: flow diagram} depicts the flow diagram of the simulation. As the simulation advances, these events occur in a sequential manner, incurring associated cash flows. The outcome of the simulation is the NPV ($V$), calculated as the discounted sum of these cash flows:{\color{black}
\begin{equation}
    V = \sum_{t=1}^{t_{\text{sim}}} CF_t e^{-r_f t \Delta T},
\end{equation}}
where $CF_t$ is the cash flow at time step $t$, {\color{black}$r_f$ is the continuously compounded risk free rate}, $\Delta T$ is the time step size in years, and $t_{\text{sim}}$ is the total number of simulation time steps. {\color{black}Since the simulation is associated with stochastic events, the resultant NPV varies with each run, even under identical parameter settings.} The rest of this section provides a detailed breakdown of the components constituting the simulation.

\subsubsection{Start \& End}
The simulation initiates with the launch of a new satellite. It progresses until the {\color{black}simulation reaches} $t_{\text{sim}}$, set as the endpoint. The mechanism for monitoring the attainment of this endpoint, $t_{\text{sim}}$, is integrated within \textit{Normal Operation}.

\subsubsection{Normal Operation}
{\color{black}Normal Operation is the period when the satellite functions as intended.} Instead of being a distinct event, Normal Operation is a transitional phase between two consecutive events, maintaining a normally operational state. During this phase, the simulation checks whether it has reached its endpoint{\color{black}, $t_{\text{sim}}$}. If {\color{black}the simulation has reached $t_{\text{sim}}$, it terminates; otherwise, it continues to the next event.}

For each time step {\color{black}under Normal Operation}, the satellite consumes propellant for station-keeping maneuvers. {\color{black}The $\Delta V$ required for station-keeping maneuver over each time step ($\Delta V_{\text{stk}}$) is defined as follows:
\begin{equation}\label{eqn: DeltaVStationKeeping}
    \Delta V_{\text{stk}}=\Delta V_{\text{stk,yr}} \Delta T, 
\end{equation}
where} $\Delta V_{\text{stk},\text{yr}}$ is $\Delta V$ required for station-keeping maneuver annually. The propellant consumption for this maneuver is determined by the rocket equation, which represents the relation between the required $\Delta V$, the performance of the propulsion subsystem, and the propellant consumption:
\begin{equation}
    m_p = m_0  \left(1-\exp\left(-\frac{\Delta V}{g_0 I_{\text{sp}}}\right)\right)
    \label{eqn: rocket equation}
\end{equation}
Here, $\Delta V$ is the $\Delta V$ value required for a specific maneuver, $m_p$ is the corresponding propellant consumption, and $m_0$ is the initial gross mass of the satellite.

As a result of {\color{black}operation, satellite missions generate revenue and costs.} The associated revenue is determined dynamically. In the simulation, the dynamic environmental factors affecting the revenue generated by satellite operations are technology obsolescence and market volatility. These factors significantly influence the decision-making process regarding the design lifetime of a system \cite{saleh2004weaving}. The revenue generated from satellite operations is calculated by multiplying these factors with the baseline revenue:{\color{black}
\begin{equation}\label{eqn: revenue}
    R\left(t\right)=R_0\Delta T\phi_{\text{mar}}\left(t\right) \phi_{\text{obs}}\left(t;t'\right),\quad t'\le t\le t_{\text{sim}},
\end{equation}}
where $t'$ denotes the reference time {\color{black}step} for the level of the applied technology. The market factor ($\phi_{\text{mar}}$) follows a geometric Brownian motion \cite{saleh2004weaving, hull2022options}. The discrete version of the stochastic differential equation is as follows:{\color{black}
\begin{equation}\label{eqn: MarketVolatility}
    \Delta \phi_{\text{mar}}=\phi_{\text{mar}}\left(t+1\right)-\phi_{\text{mar}}\left(t\right)=\mu_\text{mar}  \Delta T+\sigma_\text{mar}\epsilon\sqrt{\Delta T},\quad t\le t_{\text{sim}}.
\end{equation}
where $\mu_{\text{mar}}$ and $\sigma_{\text{mar}}$ are market drift rate and volatility, $\epsilon$ is a random variable following the standard normal distribution.} Also, the technology obsolescence factor ($\phi_{\text{obs}}$) is defined as follows \cite{saleh2004weaving}:{\color{black}
\begin{equation}
    \phi_{\text{obs}}\left(t;t'\right) = \exp\left(-\left(\frac{\left(t-t'\right)\Delta T}{\theta_{\text{obs}}}\right)^2\right),
    \label{eqn: technology obsolescence factor}
\end{equation}
where $\theta_{\text{obs}}$ controls the rate of obsolescence. Thus, revenue fluctuates according to the market factor's drift and volatility, and decays as the satellite's technology ages.}

During Normal Operation, satellites are subject to in-orbit failures caused by component malfunctions. One can model the time to a certain type of in-orbit failure using a probability distribution (e.g. Weibull distribution). A simple yet effective model to describe the overall failure behavior of a satellite is the 2-Weibull distribution model. This model can capture the dynamic characteristics of in-orbit failure types, including infant mortality and wear-out failures \cite{saleh2011spacecraft}. The reliability function ($Rel$) following the 2-Weibull distribution is expressed as follows:{\color{black}
\begin{equation}\label{eqn: WeibullMixtureReliability}
    Rel\left(t;t'\right)=\begin{cases}\alpha_{\text{rel}} \exp\left(-\left(\frac{\left(t-t'\right)\Delta T}{\theta_{\text{rel},1}}\frac{T_\text{ref}}{T_{\text{life}}}\right)^{\beta_{\text{rel},1}}\right)+\left(1-\alpha_{\text{rel}}\right) \exp\left(-\left(\frac{\left(t-t'\right)\Delta T}{\theta_{\text{rel},2}}\frac{T_\text{ref}}{T_{\text{life}}}\right)^{\beta_{\text{rel},2}}\right),& t'\le t\le t'+t_{\text{life}}\\
    0,& t > t'+t_{\text{life}}
    \end{cases}
\end{equation}}
where $T_{\text{ref}}$ is the reference design lifetime, $T_{\text{life}}$ is the design lifetime, and {\color{black}$\alpha_{\text{rel}}$, $\beta_{\text{rel},1}$, $\beta_{\text{rel},2}$, $\theta_{\text{rel},1}$, and $\theta_{\text{rel},2}$} are parameters of the 2-Weibull distribution, and {\color{black}$t'$ is the time step where the satellite begins its operation.} In this equation, the first term containing the parameter {\color{black}$\beta_{\text{rel},1}(<1)$} represents the infant mortality phase, and the second term, characterized by the parameter {\color{black}$\beta_{\text{rel},2}(>1)$}, reflects the wear-out phase. 

{\color{black}Following the reliability function, let $P_{\text{iof}}\left(t;t'\right)$ denote the probability of an in-orbit failure occurring at time step $t$ given that the satellite began operations at time step $t'$. Recognizing that the time horizon is discretized while failures can occur at any continuous point, a failure occurring between time steps $\left(t-1\right)$ and $t$ is regarded as happening at time step $t$. The probability is thus defined via the reliability function as:
\begin{equation}\label{eqn: in-orbit failure probability}
    P_{\text{iof}}\left(t;t'\right)=Rel\left(t-1;t'\right) - Rel\left(t;t'\right),\quad t'+1\le t.
\end{equation}
If a satellite that began operations at time step $t'$ is operational at time step $\left(t-1\right)$, then the conditional probability of in-orbit failure at time step $t$ is given as follows:
\begin{equation}\label{eqn: conditional in-orbit failure probability}
    P\left[t_{\text{iof}}=t;t'; t_{\text{iof}}\ge t\right]=\frac{P_{\text{iof}}\left(t;t'\right)}{Rel\left(t-1;t'\right)}=\frac{P_{\text{iof}}\left(t;t'\right)}{\sum_{\tau\ge t}P_{\text{iof}}\left(\tau;t'\right)},
\end{equation}
where $t_{\text{iof}}$ is the time step at which the in-orbit failure occurs}

\subsubsection{Decision}
\begin{figure}[hbt!]
	\centering
	\includegraphics[page=2,width=0.7\columnwidth]{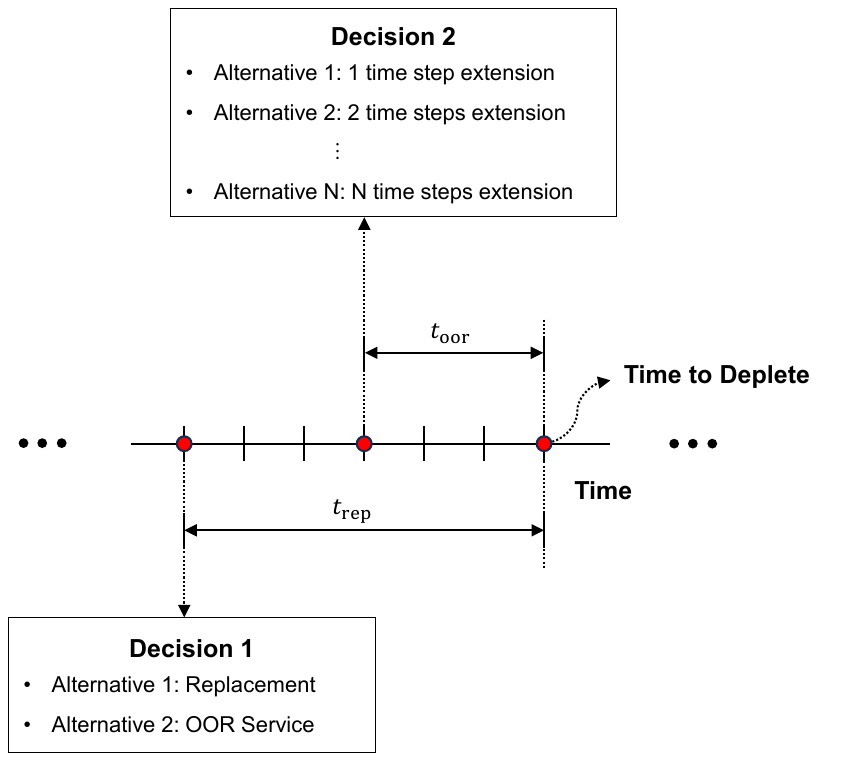}
	\caption{Decision timings and alternatives{\color{black}.}}
    \label{fig: decision timings and alternatives}
\end{figure}

Decision is a critical event that models the operational flexibility of the satellite operator. This event has two types, each with its specific timing and set of available alternatives (Fig. \ref{fig: decision timings and alternatives}). At each Decision event, the {\color{black}best} alternative is selected from the available alternatives. {\color{black}The chosen action yields varying results, depending on the interaction between the altered state of the satellite and its external environment.}

The best alternative is chosen based on the predefined utility function for each alternative. {\color{black}Assumptions regarding the satellite operator's information and preferences are established to define these utility functions:}
\begin{enumerate}
    \item The operator is fully informed about the dynamic and stochastic models of the system, including reliability, launch and service failures, market volatility, and technology obsolescence.
    \item The operator employs a greedy decision-making approach, prioritizing the financial gains before making the next decision.
\end{enumerate}
{\color{black}Based on these assumptions, the general formulation of the utility function is expressed as follows:
\begin{equation}\label{eqn: general utility}
    U\left(t\right)=\sum_{t_1=t+1}^{t_{\text{sim}}} P\left[t_{\text{E}}=t_1\right]A/P\left(r_f,t_1\right)\sum_{t_2=t+1}^{t_1}\mathbb{E}\left[CF_{t_2}\right]e^{-r_ft_2\Delta T},\quad t\le t_{\text{sim}}
\end{equation}
where $t_{\text{E}}$ is the time step where the next event occurs, $CF_{t_2}$ is the cash flow at time step $t_2$, and $A/P\left(r_f,t_1\right)$ is the equivalent annuity (EA) factor with interest rate $r_f$ over $t_1$ periods, which is defined as follows: 
\begin{equation}\label{eqn: EA continuous}
    A/P\left(r_f,t_1\right)=\frac{e^{r_ft_1} \left(e^{r_f}-1\right)}{e^{r_ft_1}-1}.
\end{equation}
Introducing the EA factor allows for fair comparison across different cash flow horizons.

Detailed formulations of the utility functions for specific decisions and their alternatives are provided in Appendices A.2 and A.3.}

Decision 1 is initiated {\color{black}$t_{\text{rep}}$} time steps before the projected propellant depletion. At that point, the operator decides whether to replace or replenish the satellite---corresponding to \textit{Replacement} and \textit{OOR Service} actions. {\color{black}If Replacement has a higher utility than OOR Service, then the Replacement is scheduled after $t_{\text{rep}}$} time steps, which is the time required to manufacture and launch for replacement.

Decision 2 is triggered after OOR {\color{black}Service} is selected in Decision 1. It occurs {\color{black}$t_{\text{oor}}$ ($<t_{\text{rep}}$)} time steps before the projected depletion. This event determines the quantity of propellant for refueling. The life-extension period provided by OOR Service is determined as the period with the highest utility, similar to Decision 1. When a specific life-extension duration is selected, OOR Service with the corresponding propellant amount ($m_{p,\text{oor}}$) is planned after {\color{black}$t_{\text{oor}}$} time steps.

\subsubsection{Action}
Action marks an event where a specific action, Replacement, or OOR Service is executed. This action is determined by the decision made in the previous Decision event or is initiated by Failure. The Action event {\color{black}is} a bridge between decision-making and operational execution, leading to a change in the {\color{black}satellite's state.}

Each action is associated with distinct failures. {\color{black}If a failure occurs during the replacement or OOR, the simulation proceeds to Failure; otherwise, the associated costs of the satellite replacement or OOR service occur. Subsequently, the system's state is updated, and the simulation progresses to the next event.}

In Replacement, the satellite currently in orbit is substituted with a new one. Replacement is successful if the launch and the orbital transfer are performed successfully. Launch is failed by a probability of {\color{black}$p_{\text{lau}}$}. In addition, the orbital transfer is successful if the initially loaded propellant is enough to accomplish the maneuvers. The following equation represents the probability of replacement failure ({\color{black}$p_{\text{rep}}$}), considering launch failure rate ({\color{black}$p_{\text{lau}}$}) and orbital transfer successes:{\color{black}
\begin{equation}
    p_{\text{rep}}=1-\left(1-p_{\text{lau}}\right) P\left[m_{p,\text{ot}}\le m_{p,\text{des}}\right],
    \label{eqn: replacement failure probability}
\end{equation}}
where $P\left[m_{p,\text{ot}} \le m_{p,\text{des}}\right]$ denotes the probability that the design propellant mass, $m_{p,\text{des}}$, is at least as large as the required propellant for the orbital transfer, $m_{p,\text{ot}}$ (i.e., the likelihood of a successful orbital transfer). To obtain this value, a stochastic model of $m_{p,\text{ot}}$ is required.

Launch service providers give the accuracy of the orbit injection to their customers. For instance, the Ariane {\color{black}6} User's Manual \cite{arianespace2021ariane} specifies the injection accuracy data of {\color{black}a} geostationary transfer orbit (GTO) mission as detailed in Table \ref{tab: injection accuracy}.
\begin{table}[hbt!]
\begin{center}
    \caption{GTO injection accuracy of Ariane {\color{black}6} \cite{arianespace2021ariane}}
    \begin{tabular}{lcc} 
        \hline\hline
        Description & Notation & Standard Deviation\\\hline
        Semi-major axis deviation (km) & $a$ & 40 \\
        Eccentricity deviation & $e$ & 4.5E-4 \\
        Inclination deviation (deg) & $i$ & 0.02 \\
        Argument of perigee deviation (deg) & $\omega$ & 0.2 \\
        Right ascension of the ascending node deviation (deg) & $\Omega$ & 0.2 \\
        \hline\hline
    \end{tabular}
    \label{tab: injection accuracy}
\end{center}
\end{table}
Correcting such errors requires additional $\Delta V$ for the orbital transfer phase. We model this as a half-normal distribution:
\begin{equation}
    \Delta V_{\text{ot}} = \Delta V_{\text{ot},\text{ideal}} + \Delta V_{\text{err}},\quad \Delta V_{\text{err}}\sim \mathcal{H}\left(\sigma_{\text{oi}}\right),
\end{equation}
where $\Delta V_{\text{ot}}$ and $\Delta V_{\text{ot},\text{ideal}}$ represent the actual and ideal $\Delta V$ required for orbital transfer from the initial orbit to the target orbit, respectively, and $\Delta V_{\text{err}}$ is the additional $\Delta V$ needed for injection error correction. Here, $\mathcal{H}\left(\sigma_{\text{oi}}\right)$ denotes a half-normal distribution derived by taking the absolute value of a normal distribution with mean $0$ and standard deviation $\sigma_{\text{oi}}$.

Combining this with Eq. \eqref{eqn: rocket equation}, the propellant mass and $\Delta V$ relationship, {\color{black}explains the second multiplier in the second term of Eq. \eqref{eqn: replacement failure probability}:
\begin{equation}
    P\left[\Delta V_{\text{ot}}\le \Delta V_{\text{des}}\right] = 2 \left(\Phi\left(\min\left(\frac{\Delta V_{\text{des}}-\Delta V_{\text{ot}}}{\sigma_{\text{oi}}}, 0\right)\right)-0.5\right),
\end{equation}
where $\Delta V_{\text{des}}=g_0 I_\text{sp} \ln \left(m_{\text{wet}}/m_{\text{dry}}\right)$.}

In {\color{black}OOR Service}, the satellite currently in orbit is refueled with the amount of propellant determined in the previous Decision 2. {\color{black}OOR Service} is subject to the risk of service failure, which occurs if any step in the series of maneuvers---including rendezvous, proximity operation, docking, undocking (RPODU)---or robotic manipulation for the servicing activity fails. These failures can lead to various levels of malfunction, ranging from minor functional degradation and delayed service to the total loss of the client satellite. In this research, we only consider catastrophic failures resulting in a total loss of the satellite, characterized by a probability of {\color{black}$p_{\text{serv}}$.}

\subsubsection{Failure}
The simulation includes three distinct types of failures to represent the realistic risks in satellite operations: in-orbit failure during {\color{black}Normal Operation}, replacement failure at {\color{black}Replacement}, and service failure {\color{black}at OOR Service}, as previously introduced. {\color{black}The probability of in-orbit failure is given by Eq. \eqref{eqn: conditional in-orbit failure probability}, the replacement failure probability by Eq. \eqref{eqn: replacement failure probability}, and the service failure rate is set to 1\% \cite{liu2021economic}. These values are summarized in Table \ref{tab: failure probability}}

\begin{table}[hbt!]
\color{black}\begin{center}
    \caption{\color{black}Failure types and corresponding probabilities}
    \begin{tabular}{lcc} 
        \hline\hline
        Description & Notation & Probability\\\hline
        In-orbit failure at $t$; givne launched at $t'$ and operational at $\left(t-1\right)$ & $P\left[t_{\text{iof}}=t;t'; t_{\text{iof}}\ge t\right]$ & Eq. \eqref{eqn: conditional in-orbit failure probability} \\
        Replacement failure & $p_{\text{rep}}$ & Eq. \eqref{eqn: replacement failure probability} \\
        Service failure & $p_{\text{serv}}$ & 0.01 \\
        \hline\hline
    \end{tabular}
    \label{tab: failure probability}
\end{center}
\end{table}

Upon any of these failures, the Replacement is planned {\color{black}$t_{\text{rep}}$} time steps after the failure occurs.

\section{Satellite System Architecting Framework}\label{sec: section 4}
The outputs from the lifecycle simulation {\color{black}serve} as substitutes for the objective functions of $\mathbf{P}_{\text{O}}$. {\color{black}Since the simulation executes a sequence of stochastic events to generate NPV outcomes,} the objective functions of $\mathbf{P}_{\text{O}}$, which consist of statistical measures of NPV, can only be estimated by repeatedly running the simulation. This sampling-based estimation, commonly known as the Monte Carlo method, presents two main difficulties in optimization: 1) the same set of design variables under the same scenario may yield different estimation results each time; and 2) conducting repetitive simulations is computationally expensive, especially given that optimization processes typically involve iterative computations.

\begin{figure}[hbt!]
	\centering
	\includegraphics[page=3,width=1.0\columnwidth]{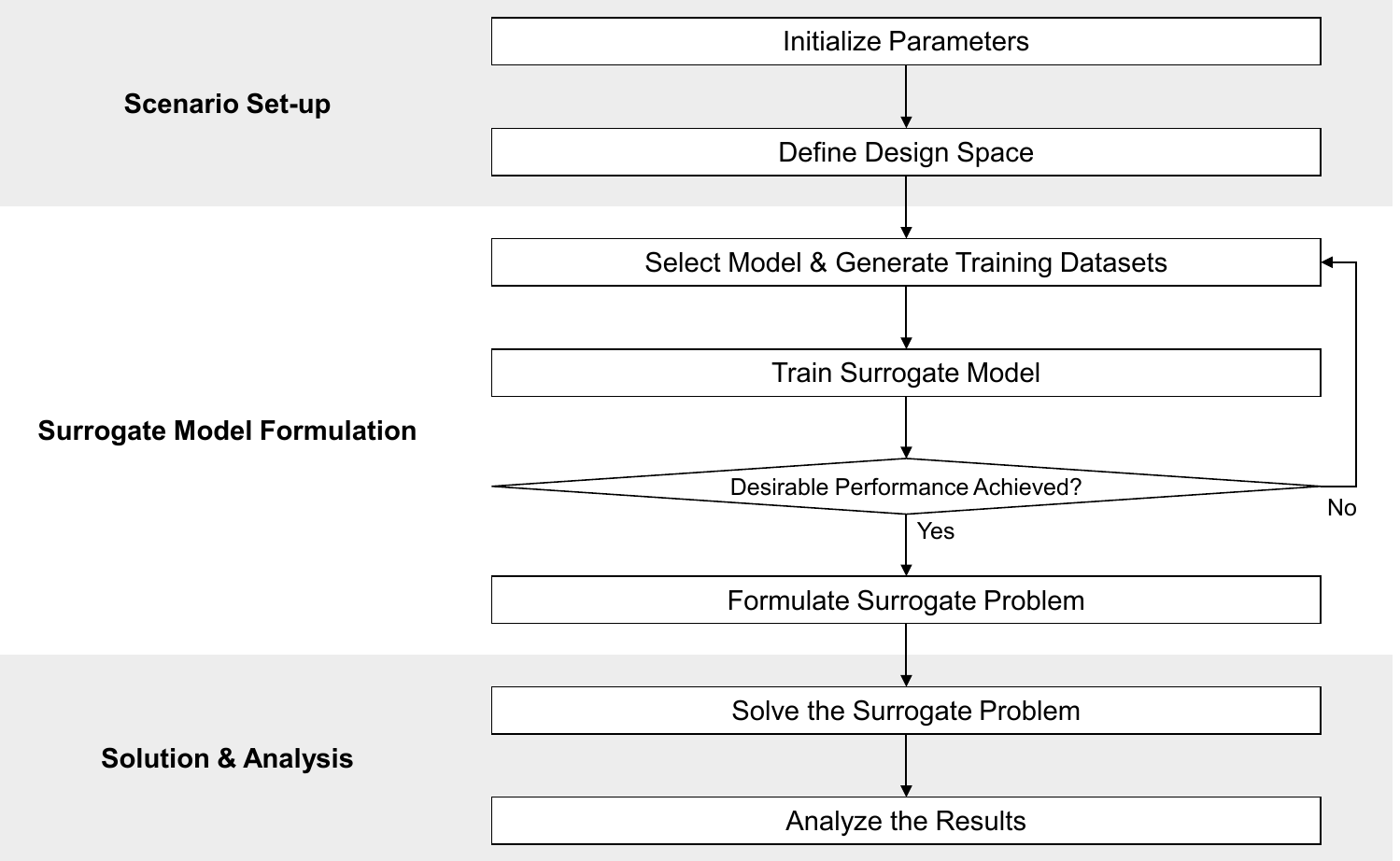}
	\caption{Proposed satellite system architecting framework{\color{black}.}}
    \label{fig: framework overview}
\end{figure}

{\color{black}To address these difficulties, we propose a simple surrogate model-based optimization framework.} Figure \ref{fig: framework overview} describes the overall procedure. The rest of this section details the processes of this framework: \textit{Scenario Set-up}, \textit{Surrogate Model Formulation}, and \textit{Solution \& Analysis}.

\subsection{Scenario Set-up}
The simulation parameters reflecting the scenario of interest are set. Following this, the design space is defined. It should encompass a meaningful range of values for the design variables while being appropriately constrained to enhance the computational efficiency of the optimization{\color{black}.}

\subsection{Surrogate Model Formulation}
After the scenario is set, the processes for generating surrogates that substitute for the objective functions are initiated. First, a {\color{black}baseline model} is selected. The alternatives range from simple linear regression with polynomials to more flexible non-parametric models, such as Gaussian process regression or neural networks. Multiple models can be selected and processed simultaneously, and the best one can be {\color{black}chosen after assessment.} For example, Gaussian processes with different kernels can be selected as candidates \cite{rasmussen2006gaussian}.

Following the model selection, a set of experiment design variables, {\color{black}$\mathcal{A}=\subset \mathcal{X}$}, is defined. This set should consist of adequate samples of the design space to ensure that the trained models can effectively mimic the simulation results. At the same time, the number of experimental points, {\color{black}$\left|\mathcal{A}\right|$}, should be regulated to balance the computational burden. Given the simulation's complexity and non-linear nature, an appropriate selection strategy should be considered \cite{montgomery2012design}.

For each element in {\color{black}$\mathcal{A}$}, $N$ Monte Carlo simulations are executed, and the sample mean, $\overline{V}_i$, and sample standard deviation, $S_i$, of NPV are used as unbiased estimations of the mean and standard deviation of the NPV:
\begin{equation}
    \overline{V}_i = \frac{1}{N}\sum_{j=1}^NV_{ij},\quad 1\le i\le {\color{black}\left| \mathcal{A}\right|},
    \label{eqn: SampleMean}
\end{equation}
\begin{equation}
    S_i = \sqrt{\frac{1}{N-1}\sum_{j=1}^N\left(V_{ij}-\overline{V}_i\right)^2},\quad 1\le i\le {\color{black}\left| \mathcal{A}\right|},
    \label{eqn: SampleStandardDeviation}
\end{equation}
where $i$ is the index for each experimental point in {\color{black}$\mathcal{A}$}, $j$ is the index for each experiment result among $N$ simulations, $\overline{V}_i$ is the sample mean of the $i$th experimental point, $S_i$ is the sample standard deviation of the $i$th experimental point, and $V_{ij}$ is the $j$th result of the $i$th experimental point. These estimations converge to corresponding values as $N$ increases, though this also escalates the computational burden, necessitating a careful determination of $N$ \cite{rice2006mathematical}.

 The results of these experiments aggregate to form datasets $\mathcal{D}_1$ and $\mathcal{D}_2$, comprising input-output pairs used for training the surrogate models for the objective functions:
\begin{equation}
    \mathcal{D}_1 = \bigl\{\bigl(\mathbf{x}_i, \overline{V}_i\bigr)\bigr\}_{1\le i\le \left| {\color{black}\mathcal{A}}\right|},\quad \mathcal{D}_2 = \bigl\{\bigl(\mathbf{x}_i, \overline{V}_i/S_i\bigr)\bigr\}_{1\le i\le \left| {\color{black}\mathcal{A}}\right|}.
\end{equation}

The performance of the trained model is assessed using the test sets. These sets are constructed similarly to the training sets but with design variables randomly sampled from the design space. The value of a performance metric, such as $R^2$, is calculated to evaluate the model's efficacy. If the performance is unsatisfactory, the process cycles back to either model selection or experiment design, repeating the procedure until satisfactory performance levels are attained.

\subsection{Solution \& Analysis}
If the surrogates for the objective functions demonstrate satisfactory performance, they are set to substitute for the objective functions of the original problem ($\mathbf{P}_{\text{O}}$):
\newline\newline\noindent
($\mathbf{P}_{\text{SO}}$) Surrogate-Assisted Optimization of Satellite System Architecture
\begin{equation}
    \underset{\mathbf{x}}{\text{maximize}}\ \mathbf{\hat{J}}\left(\mathbf{x}\right) = \left[\hat{J}_1,\hat{J}_2\right]{\color{black},}
\end{equation}
subject to
\begin{equation}
    \mathbf{x}\in {\color{black}\mathcal{X}}.
\end{equation}
In this formulation, {\color{black}$\mathcal{X}$ denotes the feasible set from {\color{black}$\mathbf{P}_{\text{O}}$}, and $\hat{J}_1$ and $\hat{J}_2$ represent the surrogates for $E_V$ and $E_V/\sigma_V$, respectively. Using} these surrogate objective functions overcomes the difficulties mentioned earlier. Firstly, they provide deterministic output; if the input value remains constant, the output value will also be consistent, thereby facilitating the optimization process. Secondly, the surrogate models yield the output through simple arithmetic processes rather than complex simulations, easing the repetitive assessment of each design choice and thus streamlining the optimization process.

{\color{black}With this formulation, }the selection of the optimization algorithm becomes relatively flexible. Metaheuristics such as NSGA-II \cite{deb2002fast} can be simple and effective to obtain the solutions.

After identifying efficient solutions, they undergo a thorough analysis. Building on the insights obtained from the analysis, a new scenario with different parameters can be established for further analysis. We showcase the detailed process and effectiveness of the framework in the next section.
{\color{black}
\section{Case Study: GEO Communication Satellite}\label{sec: section 5}
We conducted two case studies on the architecture of GEO communication satellite systems, each employing a different propulsion technology: chemical propulsion and electric propulsion. Chemical propulsion has been a reliable approach since the dawn of the space industry, offering high thrust but relatively low efficiency. In contrast, electric propulsion — an emerging technology — offers a significantly higher specific impulse, and both systems are now widely used across various missions, ranging from Earth-orbiting satellites to deep space exploration.

The remainder of this section is organized into three subsections. First, we present a case study of the baseline scenario to diagnose the current service performance in terms of system architecture and to illustrate the flow of the proposed framework. Next, we conduct parametric studies of GEO communication satellites employing chemical and electric propulsion, respectively, to identify the threshold of service performance (characterized by cost and capacity) at which a new architectural solution emerges, namely actively employing OOR service by loading less propellant at launch (propellant-reduced architecture).

\subsection{GEO Communication Satellite with Chemical Propulsion: Baseline Scenario}
The baseline scenario presents the details of the solution procedure outlined in Section IV and assesses the current target service performance from the perspective of the satellite system architecture. The process details are presented sequentially, and parameters are carefully selected to fulfill these objectives.

\subsubsection{Scenario Set-up}
Table \ref{tab: baseline parameters chemical} summarizes the system and environmental parameters for the baseline scenario. We assume the satellite employs a hydrazine monopropellant chemical propulsion system, allowing it to be refueled via Orbit Fab's proposed service. The satellite is initially launched into GTO and then maneuvers to GEO using its onboard propulsion. The GTO perigee and apogee altitudes are set to 200 km and 35,786 km, respectively. The simulation runs for 30 years with a one-week time step to ensure adequate temporal resolution.

\begin{table}[hbt!]
\begin{center}\color{black}
    \caption{\color{black}Values of parameters for the baseline scenario with chemical propulsion}
        \begin{tabular}{l c c c} 
        \hline\hline
        Parameter & Notation & Value & Unit\\\hline
        \multirow{2}{*}{Coefficients for propulsion subsystem mass} & $a_{\text{prop}}$ & 1.336 & kg$^{1/3}$\\
        & $b_{\text{prop}}$ & 0.455 & kg\\
        Servicing instrument cost & $C_{\text{serv}}$ & 0.03 & \$M \\
        Specific launch cost & $c_{\text{lau}}$ & 0.01 & \$M/kg\\
        Fix cost of OOR service & $c_{\text{oor},f}$ & 4 & \$M\\
        Variable cost coefficient of OOR service & $c_{\text{oor},v}$ & 0.16 & \$M/kg \\
        Specific impulse of the propulsion subsystem & $I_{\text{sp}}$ & 230 & sec\\
        OOR service capacity & $M_{\text{oor}}$ & 100 & kg \\
        Reference payload mass & $m_{\text{pl},\text{ref}}$ & 500 & kg\\
        Reference mass for other subsystems & $m_{\text{base},\text{ref}}$ & 600 & kg\\
        Servicing instrument mass & $m_{\text{serv}}$ & 1.1 & kg \\
        Launch failure rate & $p_{\text{lau}}$ & 0.03 & -\\
        OOR service failure rate & $p_{\text{oor}}$ & 0.01 & - \\
        Initial operational revenue & $R_0$ & 70 & \$M/yr \\
        Risk free rate & $r_f$ & $\ln\left(1+0.03\right)$ & - \\
        Time to OOR service & $T_{\text{oor}}$ & 4/52 & yr\\
        Reference design lifetime & $T_{\text{ref}}$ & 15 & yr\\
        Time to replacement & $T_{\text{rep}}$ & 3 & yr\\
        Simulation time & $T_\text{sim}$ & 30 & yr\\
        Time step & $\Delta T$ & 1/52 & yr\\
        $\Delta V$ required for orbital transfer from the ideal initial orbit to the target orbit & $\Delta V_\text{ot}$& 1477 & m/s\\
        Annual $\Delta V$ required for station-keeping & $\Delta V_{\text{stk},\text{yr}}$ & 50 & m/s\\
        ADCS mass ratio & $\alpha_{\text{adcs}}$ & 0.06 & -\\
        Insurance cost ratio & $\alpha_{\text{ins}}$ & 0.2 & - \\    
        Operational cost ratio & $\alpha_{\text{op}}$ & 0.1 & - \\
        Structural mass ratio & $\alpha_{\text{str}}$ & 0.21 & -\\
        \multirow{5}{*}{Parameters for reliability function} & $\alpha_{\text{rel}}$ & 0.9490 & - \\
        & $\beta_{\text{rel},1}$ & 0.4458 & - \\
        & $\beta_{\text{rel},2}$ & 4.6687 & - \\
        & $\theta_{\text{rel},1}$ & 39830.5 & yr \\
        & $\theta_{\text{rel},2}$ & 9.8 & yr \\
        Parameter for technology obsolescence & $\theta_{\text{obs}}$ & 20 & yr\\
        Mass growth rate of a satellite to the design lifetime & $\kappa$ & 0.03 & - \\
        Market drift rate & $\mu_{\text{mar}}$ & 0.03 & yr$^{-1}$\\
        Market volatility & $\sigma_{\text{mar}}$ & 0.1 & yr$^{-1/2}$\\
        Parameter of $\Delta V$ distribution for correcting the orbit injection error & $\sigma_{\text{oi}}$ & 25 & m/s\\
        \hline\hline
        \end{tabular}
    \label{tab: baseline parameters chemical}
\end{center}
\end{table}

The parameters for the MER—$\alpha_{\text{str}}$, $\alpha_{\text{adcs}}$, and $\kappa$—are taken from Saleh et al. \cite{saleh2002spacecraft}. Additionally, $a_{\text{prop}}$ and $b_{\text{prop}}$ are determined by fitting the model presented in Eq. \eqref{eqn: MER4} to the mass data of high Earth orbit satellites \cite{wertz2011space}. The reference masses of payload ($m_{\text{pl},\text{ref}}$) and other subsystems ($m_{\text{base},\text{ref}}$) are set to 500 kg and 600 kg, respectively. With a design lifetime of $T_{\text{life}} = 15$ years and a design propellant mass ($m_{p,\text{des}}$) 3,500 kg, the resulting dry mass ($m_{\text{dry}}$) is 1,930 kg, which is within a reasonable range. The parameters for the reliability function---$\alpha_{\text{rel}}$, $\beta_{\text{rel},1}$, $\beta_{\text{rel},2}$, $\theta_{\text{rel},1}$, and $\theta_{\text{rel},2}$---are adopted from Saleh and Castet \cite{saleh2011spacecraft}. Finally, the launch cost is set within a reasonable range, referring to current and near-future vehicles capable of delivering a GEO communication satellite to GTO---Falcon 9, Ariane 6, and Atlas V \cite{arianespace2021ariane, wall2024europe, satnews2024, atlasv, falcon9}, and the remaining parameters are set based on data from existing articles \cite{saleh2004weaving, liu2021economic, wertz2011space} to ensure they fall within a reasonable range. 

The baseline service cost and capacity are based on currently available data, including Orbit Fab’s target performance metrics \cite{orbitfab} and the progression model of Showalter et al. \cite{showalter2025progression}. The OOR failure rate is assumed to be 0.01, in line with proposed by Liu et al. \cite{liu2021economic}. Service instrument mass and cost are set with reference to Orbit Fab's refueling port specifications \cite{orbitfab}. Finally, we conservatively estimate the OOR service duration at one month, allowing sufficient time for GEO rendezvous and RPODU maneuvers and slightly exceeding estimates from existing responsiveness studies \cite{dujonchay2017quantification, ho2020semi}.

The design space $X$ is defined as Table \ref{tab: design space chemical}. Here, we select the range of the design lifetime following the range studied in \cite{saleh2002spacecraft}, which the MER and CER are based on. Additionally, the range of design propellant mass is selected to cover the minimum to the maximum, considering the given performance of the propulsion subsystem and the mass parameters.
\begin{table}[hbt!]
\begin{center}\color{black}
    \caption{\color{black}Design space for case study of chemical propulsion}
        \begin{tabular}{l c c c} 
        \hline\hline
        Design variable & Notation & Range & Unit\\\hline
        Design lifetime & $T_{\text{life}}$ & $\left[5, 15\right]$ & yr\\
        Design propellant mass & $m_{p,\text{des}}$ & $\left[1500, 3500\right]$ & kg\\
        \hline\hline
        \end{tabular}
    \label{tab: design space chemical}
\end{center}
\end{table}

\subsubsection{Surrogate Problem Formulation}
We select Gaussian processes with various kernels as baseline models for surrogates. These kernels include the squared exponential (SE), Mat\'ern 5/2, and Mat\'ern 3/2 \cite{rasmussen2006gaussian}. Each kernel corresponds to differentiability classes $C^\infty$, $C^2$, and $C^1$, respectively. The SE kernel is known for its smoothness, making it ideal for capturing continuous trends, whereas the Mat\'ern kernels are selected for their ability to model less smooth responses.

We employ a full factorial design approach to sample experimental points for training these models. Specifically, the design lifetime is segmented into 1-year intervals, and the design propellant mass is segmented into 100 kg increments. This segmentation yields 11 distinct levels for the design lifetime and 21 levels for the design propellant mass, resulting in 231 variable combinations. These combinations form the set of experimental points, $\mathcal{A}$.

For each experimental point, 400 Monte Carlo simulations are executed. As a result, the datasets $\mathcal{D}_1$ and $\mathcal{D}_2$ are generated. These datasets train all the candidate models. As a result of the training process, the hyperparameters of each model are updated to maximize the marginal likelihood \cite{rasmussen2006gaussian}.

\begin{table}[hbt!]
\begin{center}\color{black}
    \caption{\color{black}Performance of surrogates}
    \begin{tabular}{lcr} 
    \hline\hline
     
    Kernel                                  &  Objective Function   & $R^2$\\\hline
    \multirow{2.5}{*}{SE}                   & $J_1$                 & 0.940\\
                                            & $J_2$                 & 0.958\\
    \multirow{2.5}{*}{Mat\'{e}rn 5/2}       & $J_1$                 & 0.946\\
                                            & $J_2$                 & 0.961\\
    \multirow{2.5}{*}{Mat\'{e}rn 3/2}       & $J_1$                 & 0.949\\
                                            & $J_2$                 & 0.963\\
    \hline\hline
    \end{tabular}
    \label{tab: surrogates performance}
\end{center}
\end{table}

The test set assesses the trained models. The design variables for the test set are sampled uniformly at random within the design space. The performance of each model is provided in Table \ref{tab: surrogates performance}, as $R^2$. The result reveals that all kernel alternatives for each objective function achieve $R^2$ values exceeding 0.9 in the test. This result implies that the sampling strategy and baseline models are appropriate. Among the kernel alternatives, Mat\'{e}rn 3/2 exhibits the highest score, the roughest kernel among the alternatives. Consequently, the Gaussian process with Mat\'{e}rn 3/2 is selected for the surrogate model to be used in the subsequent optimization process.

Figure \ref{fig: baseline heatmap} provides contour plots with heatmaps for the selected surrogate model corresponding to each objective function.

\begin{figure}[hbt!]
	\centering
	\includegraphics[page=4,width=1.0\columnwidth]{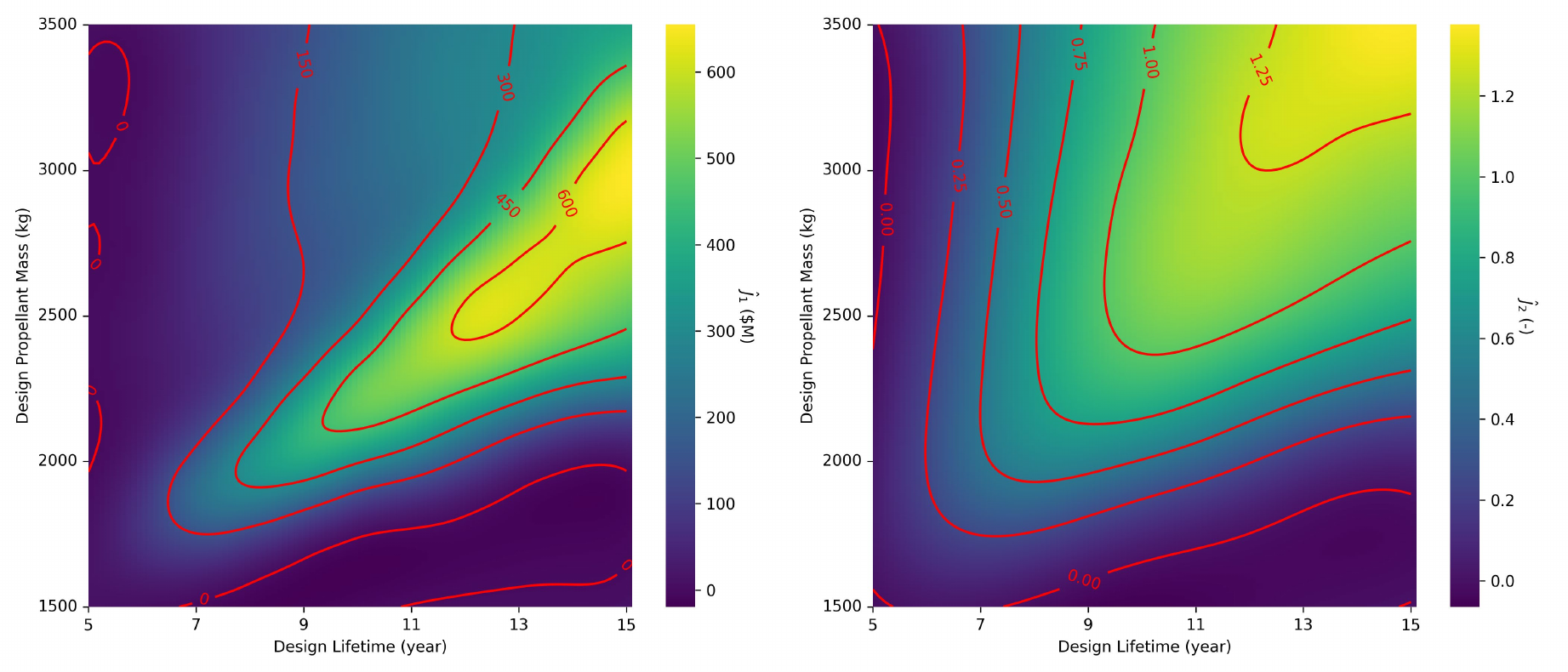}
	\caption{\color{black}Contour plot and heatmap of $\hat{J}_1$ (left) and $\hat{J}_2$ (right) for the baseline scenario.}
    \label{fig: baseline heatmap}
\end{figure}

\subsubsection{Solution \& Analysis}
A set of efficient solutions is obtained using NSGA-II \cite{deb2002fast} as the optimization algorithm. Figure \ref{fig: baseline result plot} depicts these solutions graphically in design and normalized objective spaces. In the objective space, the utopia point is marked by a red star, indicating the ideal values of the objective functions.

\begin{figure}[hbt!]
	\centering
	\includegraphics[page=5,width=1.0\columnwidth]{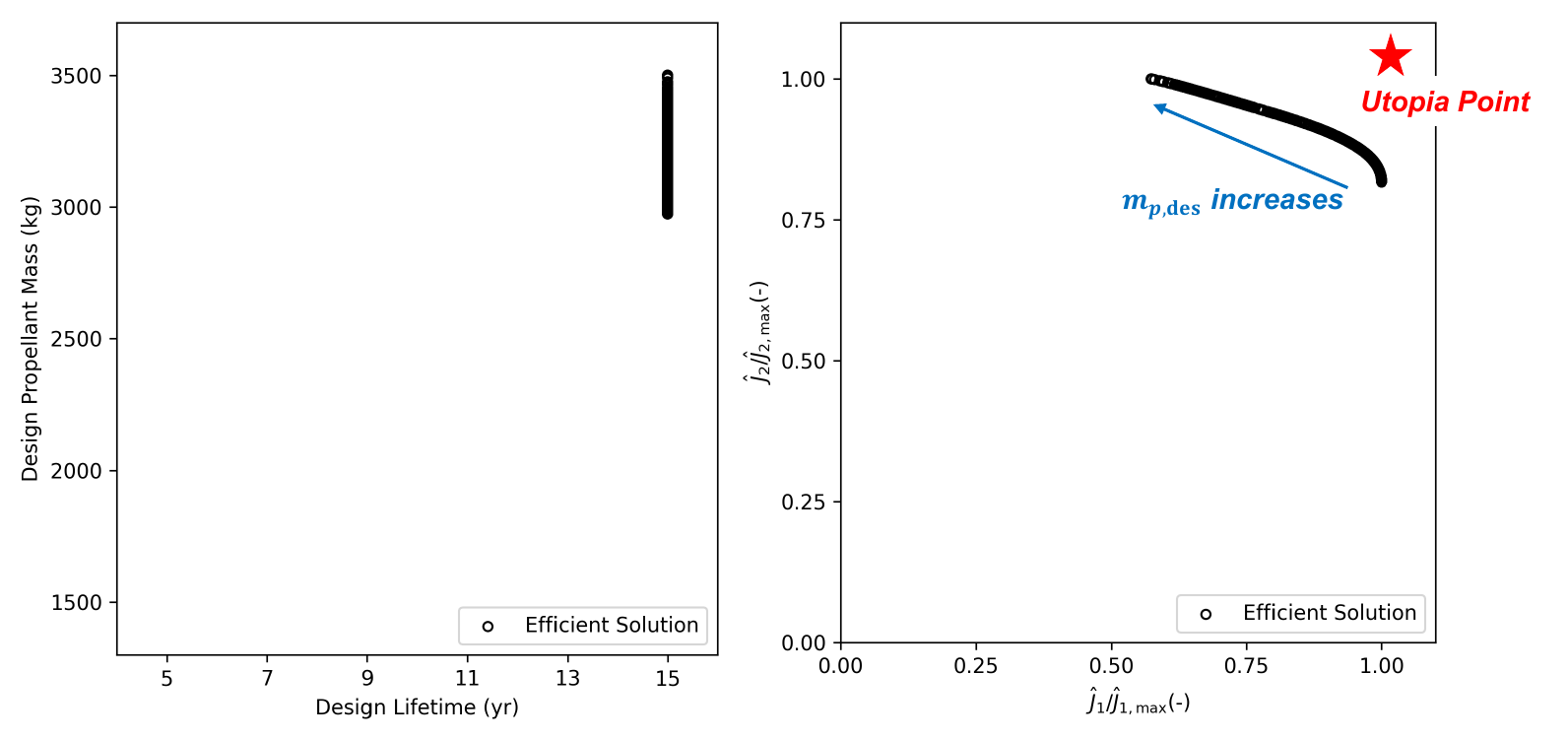}
	\caption{\color{black}Pareto front of baseline scenario in design space (left) and normalized objective space (right).}
    \label{fig: baseline result plot}
\end{figure}

In this baseline scenario, all solutions converge on a 15-year design lifetime, with propellant mass varying only between about 3,000 kg and 3,500 kg. Although a trade-off exists---loading more propellant at launch reduces overall profitability but raises profit per unit of risk---every design retains sufficient propellant throughout its lifetime.

This convergence reveals a unified design strategy: maximize lifetime while ensuring enough propellant reserves. Consequently, all solutions form a single architectural class, in line with current satellite design and operational trends. Under the given service parameters, further architectural modifications thus offer minimal additional benefit.

\subsection{GEO Communication Satellite with Chemical Propulsion: Parametric Study}
In the baseline scenario, using the current target service parameters, only the traditional design approach remains efficient. Since OOR service is still emerging in the industry and is at an early stage, its performance is expected to improve over time \cite{showalter2025progression}. But when will a new architectural solution emerge? Equivalently, at what service performance threshold should a satellite operator consider reducing on-board propellant mass at launch and actively sourcing propellant in orbit? To answer these questions, we conduct a parametric study varying two key parameters: service capacity and cost.

We repeat the procedure performed in the baseline scenario for every combination of service capacity and cost as presented in Table \ref{tab: chemical parametric study service parameters}. Capacity is scaled from the baseline level (index 1) to 4 (index 2), 7 (index 3), and 10 (index 4). Service cost is scaled from the baseline (index 1) down to 0.8, 0.6, 0.4, and 0.2 (indices 2–5). Pairing each of the four capacity levels with each of the five cost levels yields a total of 20 scenarios.

\begin{table}[hbt!]
\begin{center}\color{black}
    \caption{\color{black}Service cost and capacity parameters for parametric study of chemical propulsion by scenario index}
    \begin{tabular}{c c c c}
    \hline\hline
    Index & Service Capacity & \multicolumn{2}{c}{Service Cost} \\\cline{2-4}
    & $M_{\text{oor}}$, kg & $c_{\text{oor},f}$, \$M & $c_{\text{oor},v}$, \$M/kg \\\hline
    1 & 100 & 4 & 0.16 \\
    2 & 400 & 3.2 & 0.128 \\
    3 & 700 & 2.4 & 0.096 \\
    4 & 1,000 & 1.6 & 0.064 \\
    5 & -- -- & 0.8 & 0.032 \\
    \hline\hline
    \end{tabular}
    \label{tab: chemical parametric study service parameters}
\end{center}
\end{table}

Figure \ref{fig: chemical parametric study result} summarizes the outcomes. Black dots indicate scenarios in which only the conventional architecture remains efficient, while red stars mark the two cases---(capacity index 3, cost index 5) and (capacity index 4, cost index 5)---where reducing initial propellant mass becomes efficient. Based on these results, we draw the service performance threshold (blue dotted line) that separates the conventional regime from the emergent architecture.

\begin{figure}[hbt!]
	\centering
	\includegraphics[page=6,width=0.8\columnwidth]{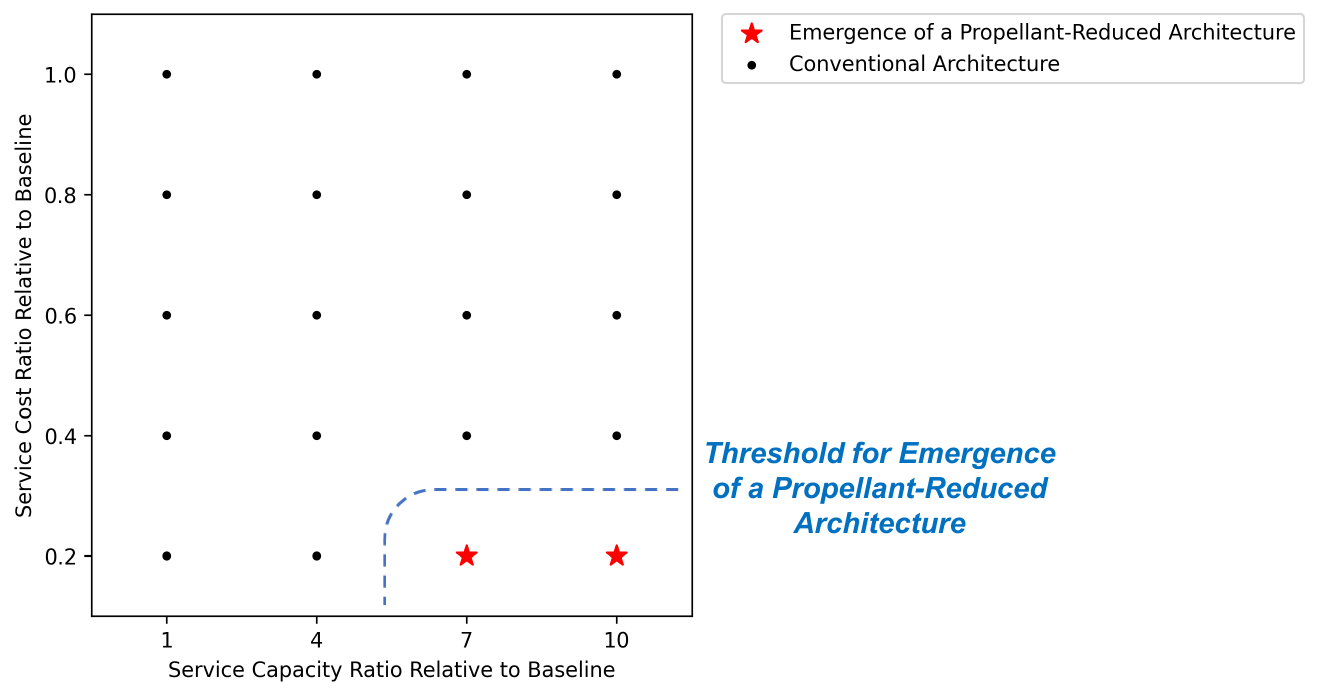}
	\caption{\color{black}Results of parametric study of GEO communication satellite with chemical propulsion.}
    \label{fig: chemical parametric study result}
\end{figure}

Figures \ref{fig: chemical parametric study result: 3-5} and \ref{fig: chemical parametric study result: 4-5} show detailed results for the two cases in which a reduced-propellant architecture emerges. In both cases, an architecture with a significantly lower propellant mass than the conventional design appears in the bottom-right corner of the objective space, indicating that this solution is both the most profitable and the riskiest. In this architecture, the spacecraft is launched with only a partial amount of propellant for its design lifetime and is subsequently replenished by OOR service. This approach becomes increasingly profitable when service performance reaches a capacity above 400 kg and cost falls below 40\% of the current target price. However, these additional servicing activities introduce further risk due to potential service failures; therefore, this solution entails the highest risk.

\begin{figure}[hbt!]
	\centering
	\includegraphics[page=7,width=1.0\columnwidth]{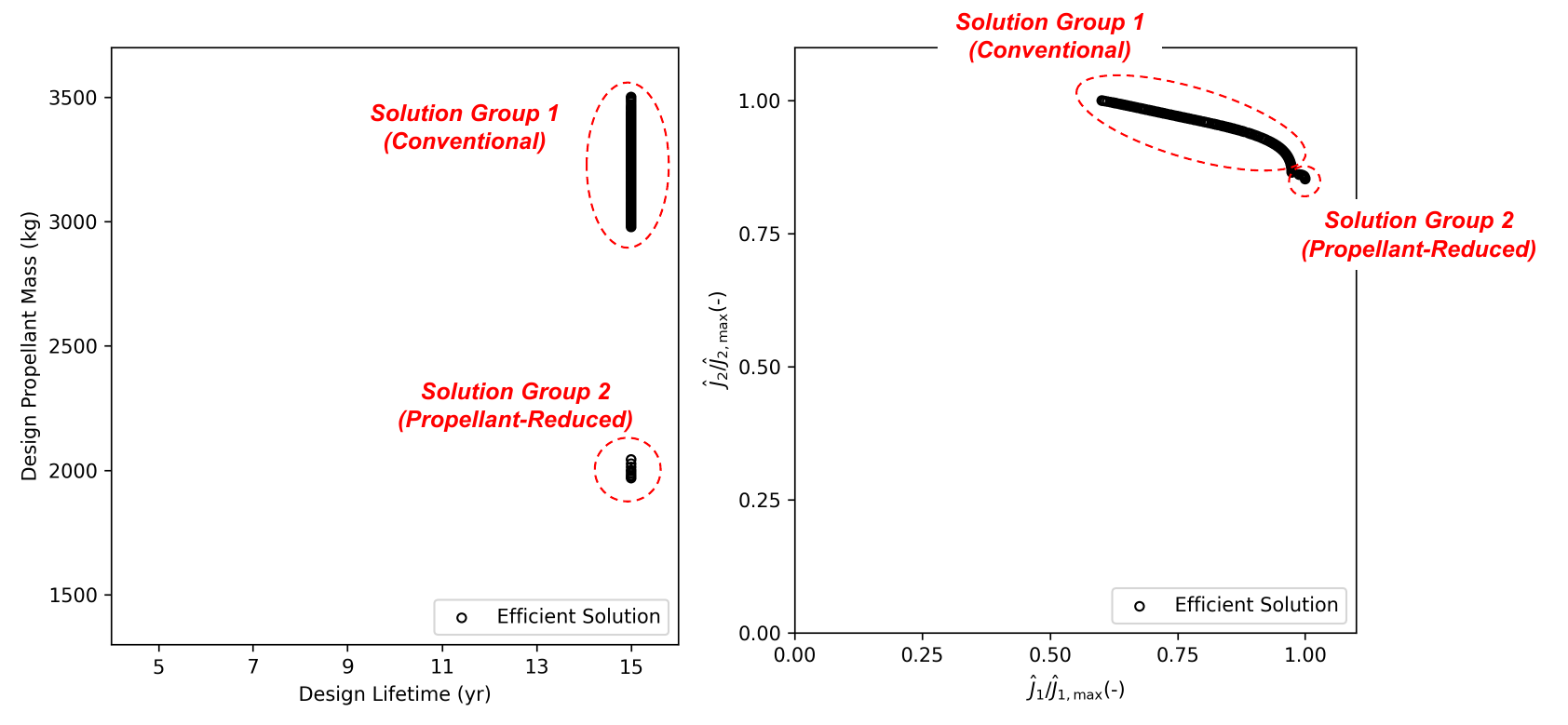}
	\caption{\color{black}Pareto front in design space (left) and normalized objective space (right): parametric study with chemical propulsion (capacity index 3 \& cost index 5).}
    \label{fig: chemical parametric study result: 3-5}
\end{figure}

\begin{figure}[hbt!]
	\centering
	\includegraphics[page=8,width=1.0\columnwidth]{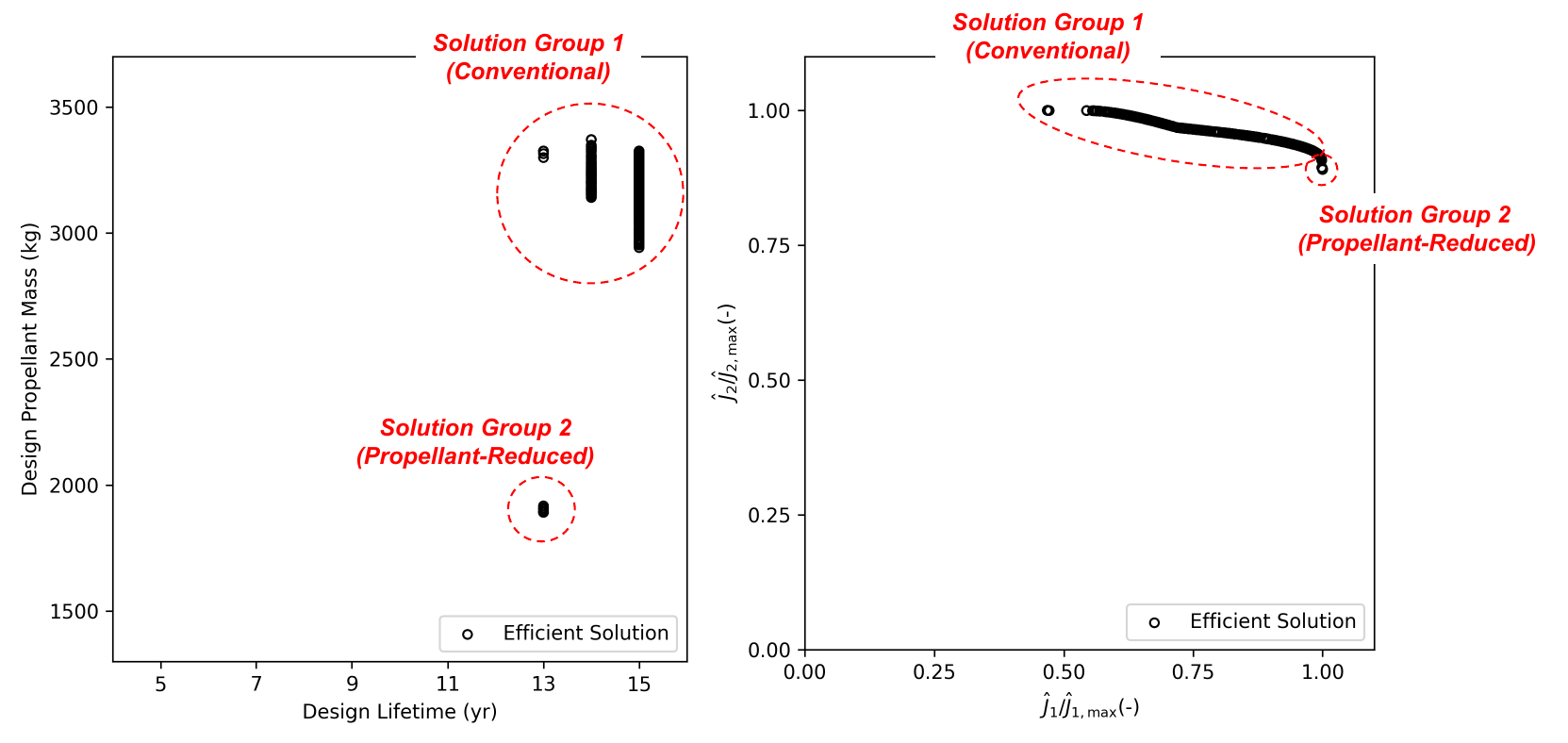}
	\caption{\color{black}Pareto front in design space (left) and normalized objective space (right): parametric study with chemical propulsion (capacity index 4 \& cost index 5).}
    \label{fig: chemical parametric study result: 4-5}
\end{figure}

\subsection{GEO Communication Satellite with Electric Propulsion: Parametric Study}
In addition to chemical propulsion, we conducted a case study of a GEO communications satellite equipped with electric propulsion---an increasingly popular option due to its high efficiency. Electric thrusters achieve specific impulses ranging from approximately 1,000 s to over 3,000 s, depending on the technology employed \cite{esa2004spt100, goebel2009evaluation}. Consequently, for a given propellant mass, electric propulsion delivers significantly greater $\Delta V$ than chemical systems. As a result, a GEO communications satellite using electric propulsion requires substantially less propellant for all mission maneuvers, from initial orbit acquisition through station-keeping.

Because electric propulsion is far more efficient, the design propellant mass influences overall system sizing in a fundamentally different way. Specifically, propellant mass represents a much smaller fraction of the total system mass; consequently, reducing the initial propellant load has only a minimal effect on sizing, and improvements in service performance yield correspondingly smaller impacts on the system architecture. To investigate this, we define the baseline scenario parameters in Table \ref{tab: baseline parameters electric}.

We assume the satellite is launched directly to GEO and carries a communication payload slightly larger than that of the chemical propulsion baseline. Accordingly, we adjust the launch cost ($c_{\mathrm{lau}}$), specific impulse ($I_{\mathrm{sp}}$), reference masses ($m_{\mathrm{pl,ref}}$, $m_{\mathrm{base,ref}}$), baseline revenue ($R_0$), and required $\Delta V$ for orbital transfer ($V_{\mathrm{ot}}$). Also, we define the design space as Table \ref{tab: design space electric}.

\begin{table}[hbt!]
\begin{center}\color{black}
    \caption{\color{black}Values of parameters for the baseline scenario with electric propulsion}
        \begin{tabular}{l c c c} 
        \hline\hline
        Parameter & Notation & Value & Unit\\\hline
        \multirow{2}{*}{Coefficients for propulsion subsystem mass} & $a_{\text{prop}}$ & 1.336 & kg$^{1/3}$ \\
        & $b_{\text{prop}}$ & 0.455 & kg \\
        Servicing instrument cost & $C_{\text{serv}}$ & 0.03 & \$M \\
        Specific launch cost & $c_{\text{lau}}$ & 0.025 & \$M/kg \\
        Fix cost of OOR service & $c_{\text{oor},f}$ & 4 & \$M \\
        Variable cost coefficient of OOR service & $c_{\text{oor},v}$ & 0.16 & \$M/kg \\
        Specific impulse of the propulsion subsystem & $I_{\text{sp}}$ & 3400 & sec\\
        OOR service capacity & $M_{\text{oor}}$ & 100 & kg \\
        Reference payload mass & $m_{\text{pl},\text{ref}}$ & 600 & kg\\
        Reference mass for other subsystems & $m_{\text{base},\text{ref}}$ & 900 & kg\\
        Servicing instrument mass & $m_{\text{serv}}$ & 1.1 & kg \\
        Launch failure rate & $p_{\text{lau}}$ & 0.03 & -\\
        OOR service failure rate & $p_{\text{oor}}$ & 0.01 & - \\
        Initial operational revenue & $R_0$ & 80 & \$M/yr \\
        Risk free rate & $r_f$ & $\ln\left(1+0.03\right)$ & - \\
        Time to OOR service & $T_{\text{oor}}$ & 4/52 & yr\\
        Reference design lifetime & $T_{\text{ref}}$ & 15 & yr\\
        Time to replacement & $T_{\text{rep}}$ & 3 & yr\\
        Simulation time & $T_\text{sim}$ & 30 & yr\\
        Time step & $\Delta T$ & 1/52 & yr\\
        $\Delta V$ required for orbital transfer from the ideal initial orbit to the target orbit & $\Delta V_\text{ot}$& 0 & m/s\\
        Annual $\Delta V$ required for station-keeping & $\Delta V_{\text{stk},\text{yr}}$ & 50 & m/s\\
        ADCS mass ratio & $\alpha_{\text{adcs}}$ & 0.06 & -\\
        Insurance cost ratio & $\alpha_{\text{ins}}$ & 0.2 & - \\    
        Operational cost ratio & $\alpha_{\text{op}}$ & 0.1 & - \\
        Structural mass ratio & $\alpha_{\text{str}}$ & 0.21 & -\\
        \multirow{5}{*}{Parameters for reliability function} & $\alpha_{\text{rel}}$ & 0.9490 & - \\
        & $\beta_{\text{rel},1}$ & 0.4458 & - \\
        & $\beta_{\text{rel},2}$ & 4.6687 & - \\
        & $\theta_{\text{rel},1}$ & 39830.5 & yr \\
        & $\theta_{\text{rel},2}$ & 9.8 & yr \\
        Parameter for technology obsolescence & $\theta_{\text{obs}}$ & 20 & yr\\
        Mass growth rate of a satellite to the design lifetime & $\kappa$ & 0.03 & - \\
        Market drift rate & $\mu_{\text{mar}}$ & 0.03 & yr$^{-1}$\\
        Market volatility & $\sigma_{\text{mar}}$ & 0.1 & yr$^{-1/2}$\\
        Parameter of $\Delta V$ distribution for correcting the orbit injection error & $\sigma_{\text{oi}}$ & 25 & m/s\\
        \hline\hline
        \end{tabular}
    \label{tab: baseline parameters electric}
\end{center}
\end{table}

\begin{table}[hbt!]
\begin{center}\color{black}
    \caption{\color{black}Design space for case study of electric propulsion}
        \begin{tabular}{l c c c} 
        \hline\hline
        Design variable & Notation & Range & Unit\\\hline
        Design lifetime & $T_{\text{life}}$ & $\left[5, 15\right]$ & yr\\
        Design propellant mass & $m_{p,\text{des}}$ & $\left[10, 50\right]$ & kg\\
        \hline\hline
        \end{tabular}
    \label{tab: design space electric}
\end{center}
\end{table}

Based on this electric propulsion baseline, we perform a parametric study by scaling the OOR service cost from 2.0 times the baseline (index 1) down to 0.2 times the baseline (index 10) in uniform decrements of 0.2---so that index 2 corresponds to 1.8 times, index 3 corresponds to 1.6 times, and so on. For each scenario, we follow the same procedure as in the chemical propulsion scenarios, but sample the training set as follows: the design lifetime is segmented into one-year intervals, and the design propellant mass is divided into 2 kg increments from 10 kg to 50 kg. As a result, these discretizations yield a total of 231 variable combinations, which form the set $\mathcal{A}$. The service cost parameters and results of each scenario are summarized in Table \ref{tab: electric parametric study service parameters}.

\begin{table}[hbt!]
\begin{center}\color{black}
    \caption{\color{black}Service cost parameters and results for parametric study of electric propulsion by scenario index}
    \begin{tabular}{c c c c}
    \hline\hline
    Index & \multicolumn{2}{c}{Service Cost} & Emergence of a Propellant-Reduced Architecture \\\cline{2-3}
    & $c_{\text{oor},f}$, \$M & $c_{\text{oor},v}$, \$M/kg & \\\hline
    1 & 8 & 0.32 & No\\
    2 & 7.2 & 0.288 & No\\
    3 & 6.4 & 0.256 & No\\
    4 & 5.6 & 0.224 & No\\
    5 & 4.8 & 0.192 & No\\
    6 & 4 & 0.16 & No\\
    7 & 3.2 & 0.128 & No\\
    8 & 2.4 & 0.096 & No\\
    9 & 1.6 & 0.064 & No\\
    10 & 0.8 & 0.032 & No\\
    \hline\hline
    \end{tabular}
    \label{tab: electric parametric study service parameters}
\end{center}
\end{table}

As expected, no new architectural solutions emerge. Figures \ref{fig: electric parametric study result: 1} and \ref{fig: electric parametric study result: 10} present detailed results for the two edge cases (indices 1 and 10), with the efficient frontier remaining unchanged. Although OOR service has minimal impact on the architecture of GEO communications satellites using electric propulsion, OOS for replacement of heavier modules---such as batteries and solar arrays---could nonetheless drive significant architectural shifts, since these subsystems constitute a much larger share of the overall system mass.

\begin{figure}[hbt!]
	\centering
	\includegraphics[page=9,width=1.0\columnwidth]{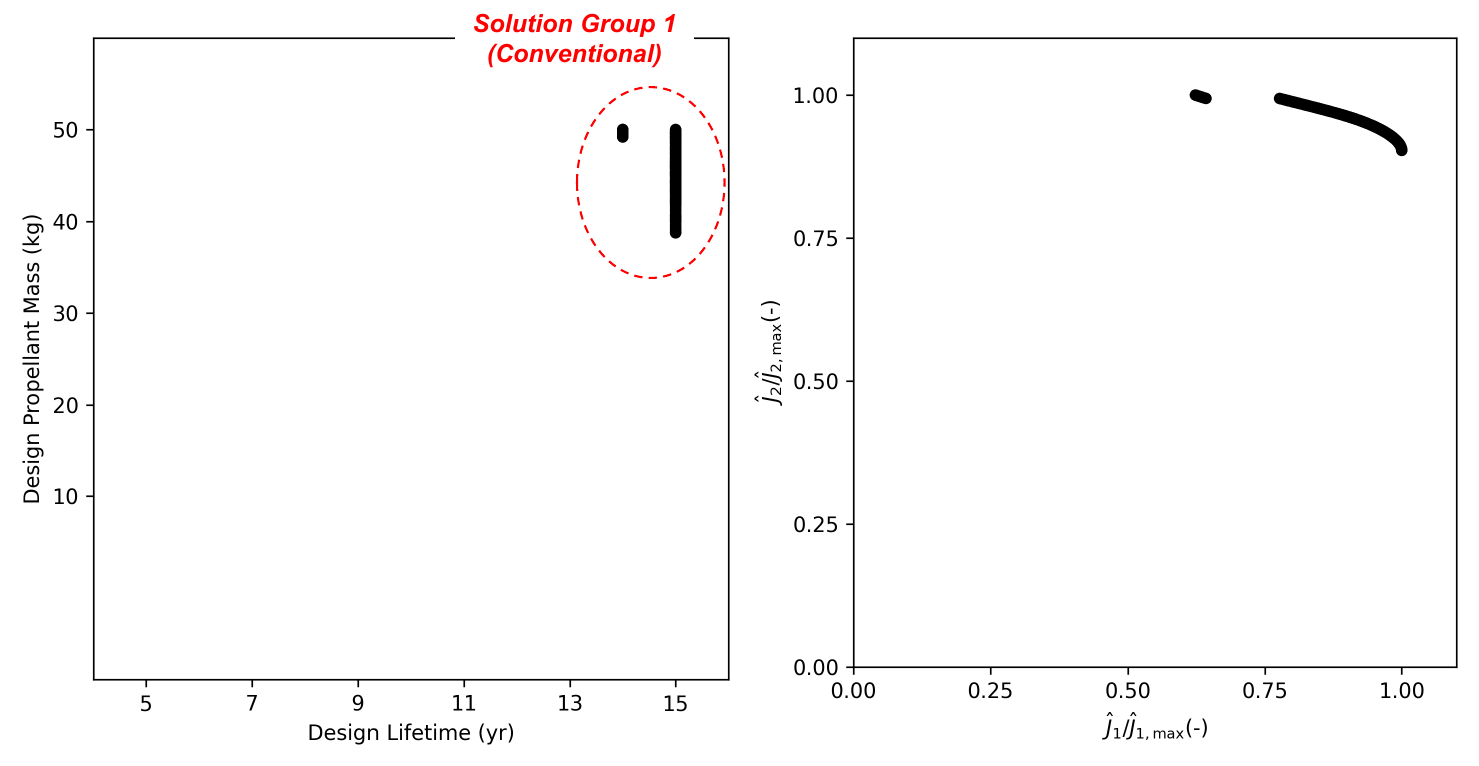}
	\caption{\color{black}Pareto front in design space (left) and normalized objective space (right): parametric study with electric propulsion (cost index 1).}
    \label{fig: electric parametric study result: 1}
\end{figure}

\begin{figure}[hbt!]
	\centering
	\includegraphics[page=10,width=1.0\columnwidth]{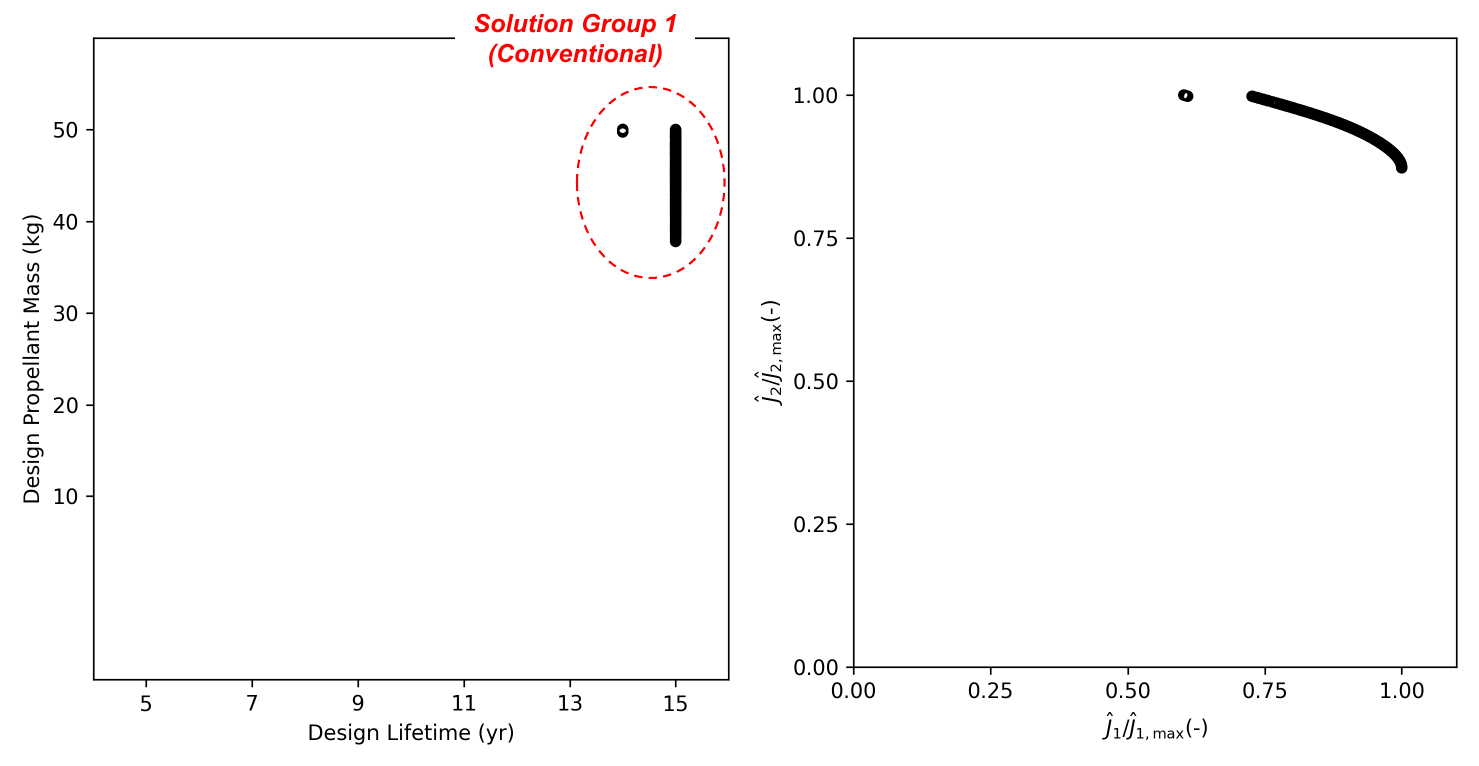}
	\caption{\color{black}Pareto front in design space (left) and normalized objective space (right): parametric study with electric propulsion (cost index 10).}
    \label{fig: electric parametric study result: 10}
\end{figure}
}

\section{Conclusions}\label{sec: section 6}
This paper introduces a satellite system architecting (high-level design) problem with explicit consideration of commercial OOR service. The design problem, formulated as a bi-objective optimization, can answer the questions of {\color{black}``What design lifetime should the satellite have?'' and ``How much propellant should be carried at launch?''} This optimization problem adopted the design lifetime and propellant as design variables, with return and risk metrics analogous to those in Portfolio Selection theory as dual objective functions. A surrogate model-based solution framework for this problem was developed based on a satellite lifecycle simulation. The simulation incorporated uncertainties and operational flexibilities, and integrated satellite sizing and cost models---revised from conventional MER and CER---to reflect the architectural implications of OOR. Furthermore, an independent servicing cost model was adopted, independent of a specific servicing infrastructure. The whole aggregate of the proposed formulation and solution methodology overcame the limitations of traditional approaches for evaluating satellite system architectures with OOR.

We conducted a case study on GEO communication satellites {\color{black} with two different propulsion systems} to demonstrate the effectiveness of our framework and . {\color{black}In this study, we first detail the solution procedure using a chemical propulsion baseline and then perform a parametric analysis for both chemical and electric configurations. Our results show that OOR service can drive significant architectural changes in GEO satellites with chemical propulsion, especially if service capacity reaches several hundred kilograms and cost drops about one-third of its current target. By contrast, for GEO satellites with electric propulsion, where propellant mass represents only a small fraction of total spacecraft mass, OOR service alone yields minimal potential for architectural change.

This study can be extended in several ways. First, by incorporating additional elements into the framework, it becomes possible to evaluate architectural solutions for a wide range of space systems under various operational environments. For example, one could model mass-constrained scenarios, where launch vehicle capacity limits the spacecraft's wet mass, by integrating appropriate wet mass constraints into the optimization problem. Although we did not examine mass-limited cases in this study (our satellite falls within the payload capacity of the launcher classes we considered), many missions face strict mass budgets, which can significantly restrict spacecraft design.

Second, this framework can be broadened to encompass various classes of OOS, including maintenance and upgrade operations. Foundational research to identify replaceable or repairable subsystems (e.g., structural and electrical components) compatible with emerging OOS technologies, perform detailed reliability-sizing analyses for each subsystem and assess performance degradation over time, and evaluate design modifications for advanced servicing (e.g., modular architectures that facilitate servicing activities) can be considered as future study subjects. By extending the framework to accommodate different OOS classes, one could derive meaningful insights into the system architecture of satellites equipped with electric propulsion systems or for smaller LEO communication satellites, considering the full potential of OOS. Ultimately, these extensions would support more informed decision-making and foster profitable, sustainable space utilization over the long term.
}

\section*{Appendix}
\subsection{Utility Function of Each Alternative}
\subsubsection{Preliminaries}
The DCF method is one of the most well-known appraisal approaches for engineering projects, classified under the income approach of business valuation \cite{newnan2017engineering}. This method involves a detailed estimation of the project's expected net cash flows---referred to as Free Cash Flow (FCF)---throughout its lifecycle. Given that FCF is projected over a specific time horizon, it is essential to account for the time value of money, as the value of identical cash differs between the present and the future. To address this, each projected cash flow is discounted to its present value. Moreover, due to the inherent uncertainty in FCF projections, the DCF method also incorporates risk considerations alongside the time value of money. Consequently, FCF is converted to its present value using a discount rate that captures both the time value of money and the associated risks, thereby yielding an estimate of the project's present economic value. Since we model associated risks in the events of the simulation, {\color{black}we use a continuously compounded risk-free rate, $r_f$, that does not consider the risk, and the NPV is represented as
\begin{equation}\label{eqn: NPV continuous}
    NPV = \sum_{T\in \mathcal{T}} CF_T e^{-r_fT},
\end{equation}
where $NPV$ is the present value of the project, $CF_T$ is the free cash flow at time $T$, and $\mathcal{T}$ is the set of times at which cash flows occur.}

Let $R\left(t_0\right)$ denote the revenue provided by an operational satellite at time step $t_0$. The expected revenue at time step $t$ is represented as{\color{black}
\begin{equation}
    \mathbb{E}\left[R\left(t\right);t', t_0\right]=R\left(t_0\right) \mathbb{E}\left[\phi_{\text{mar}}\left(t-t_0\right)\right] \phi_{\text{obs}}\left(t;t'\right),\quad t',t_0\le t,
\end{equation}
where $t'$ denotes the time step at which the satellite of interest begins operating.}

The expected {\color{black}profit} the operational satellite generates from Decision 1 until the next {\color{black}Action event, either Replacement or OOR Service}, is denoted by $EP_1$. Its mathematical expression is given by{\color{black}
\begin{equation}\label{eqn: ER1}
    EP_{1}\left(t;t'\right)=\sum_{t_1=1}^{t_{\text{rep}}}\left(\left(1-r_{\text{op}}\right)\mathbb{E}\left[R\left(t+t_1\right);t',t\right] \frac{\sum_{\tau>t+t_1}P_{\text{iof}}\left(\tau;t'\right)}{\sum_{\tau>t}P_{\text{iof}}\left(\tau;t'\right)} e^{-r_f t_1\Delta T}\right),\quad t'\le t,
\end{equation}
where $\left(1-r_{\text{op}}\right)\mathbb{E}\left[R\left(t+t_1\right);t',t\right]$ denotes the expected cash flow at time step $t+t_1$, conditional on the operational satellite not having failed; $\frac{\sum_{\tau>t+t_1}P_{\text{iof}}\left(\tau;t'\right)}{\sum_{\tau>t}P_{\text{iof}}\left(\tau;t'\right)}$ represents the conditional probability that the satellite remains operational at time step $t+t_1$ given that it is operational at time step $t$; and $e^{-r_ft_1\Delta T}$ is the continuously compounded discount factor from $t+t_1$ back to $t$, accounting for the time value of money.}

Similarly, the expected {\color{black}profit} the operational satellite generates from Decision 2 until the next action is denoted by $EP_2$:{\color{black}
\begin{equation}\label{eqn: ER2}
    EP_{2}\left(t;t'\right)=\sum_{t_1=1}^{t_{\text{oor}}}\left(\left(1-r_{\text{op}}\right)\mathbb{E}\left[R\left(t+t_1\right);t',t\right] \frac{\sum_{\tau>t+t_1}P_{\text{iof}}\left(\tau;t'\right)}{\sum_{\tau>t}P_{\text{iof}}\left(\tau;t'\right)} e^{-r_ft_1\Delta T}\right),\quad t'\le t.
\end{equation}}

\subsubsection{Decision 1}
In Decision 1, two alternatives are considered: replacement and OOR. The utility function for each alternative is defined according to the assumptions in Section \ref{sec: section 3}. We define the utility function of the Replacement ($U_{\text{rep}}$) as the expected EA by the replacement. This utility function is defined as the sum of three distinct components, reflecting various possible outcomes of the replacement:{\color{black}
\begin{align}
    U_{\text{rep};t'}\left(t\right)=U_{\text{rep,1}}\left(t;t'\right)+U_{\text{rep,2}}\left(t;t'\right),
    \label{eqn: UtilityOfReplacement}
\end{align}}
where $U_{\text{rep},1}$ represents the EA of the scenario where the replacement fails; {\color{black}$U_{\text{rep},2}$ represents the EA of the situation where the replacement is successful. The mathematical representation of each term is as follows:
\begin{equation}\label{eqn: UtilityOfReplacement1}
    U_{\text{rep,1}}\left(t;t'\right)=p_{\text{rep}} A/P\left(r_f,t_{\text{rep}}\right) EP_{1}\left(t;t'\right),
\end{equation}
\begin{align}
    U_{\text{rep,2}}\left(t;t'\right)=\left(1-p_{\text{lau}}\right)\sum_{t_1=t_{\text{rep}}+1}^{t_{\text{rep}}+t_{\text{life}}+1}P_{\text{E},1}\left(t_1;t_{\text{rep}}\right)A/P\left(r_f,t_1\right)&\Biggl(EP_{1}\left(t;t'\right)-C_{\text{IOC}}e^{-r_f t_{\text{rep}}\Delta T}\nonumber\\
    &+\sum_{t_2=t_{\text{rep}}+1}^{t_1-1}\mathbb{E}\left[R\left(t+t_2\right);t_{\text{rep}},t\right] e^{-r_ft_2\Delta T}\Biggr)\label{eqn: UtilityOfReplacement2},
\end{align}
where $P_{\text{E}}\left(t;t'\right)$ is the probability that the next event occurs at time step $t$ given $t'$, defined as follows:
\begin{align}
    P_{\text{E},1}\left(t;t'\right)=&P\left[\Delta V_{\text{ot}}+\left(t-t'\right)\Delta V_{\text{stk}}\le \Delta V_{\text{des}}\right]P_{\text{iof}}\left(t;t'\right) \nonumber\\
    &+P\left[\Delta V_{\text{ot}}+\left(t-t'-1\right)\Delta V_{\text{stk}} \le \Delta V_{\text{des}} < \Delta V_{\text{ot}}+\left(t-t'\right)\Delta V_{\text{stk}}\right]\sum_{\tau>t}P_{\text{iof}}\left(\tau;t'\right) \nonumber\\
    &-P\left[\Delta V_{\text{ot}}+\left(t-t'-1\right)\Delta V_{\text{stk}} \le \Delta V_{\text{des}} < \Delta V_{\text{ot}}+\left(t-t'\right)\Delta V_{\text{stk}}\right]P_{\text{iof}}\left(t;t'\right). \label{eqn: P_E}
\end{align}}

The utility function of the OOR Service is defined as the maximum of the functions {\color{black}$u_{\text{oor}}(k,t;t')$ with respect to $k$:
\begin{equation}
    U_{\text{oor}}(t;t')=\max_{k\in \mathcal{K}_t} u_{\text{oor}}\left(k,t;t'\right),
\end{equation}
where $k$ denotes the number of time steps by which refueling extends the satellite's operational lifetime, $u_{\text{oor}}(k,t;t')$ is the EA achieved by extending the satellite's operational lifetime by $k$ time steps, and $\mathcal{K}_t$ is available values of $k$ at time step $t$.}

The required mass of propellant for extending the operational lifetime by $k$ time steps is given by:{\color{black}
\begin{equation}
    m_{p,\text{oor}}\left(k,t\right) = \left( m_{\text{dry}} + m_{p,\text{rem}}\left(t\right) \right)  \exp\left(-\frac{\Delta V_{\text{stk}}  t_{\text{rep}}}{g_0 I_{\text{sp}}}\right)  \left( 1 - \exp\left(-\frac{\Delta V_{\text{stk}}  k}{g_0 I_{\text{sp}}}\right) \right),\quad k\in \mathcal{K}_t,
\end{equation}}
where $m_{p,\text{rem}}\left(t\right)$ is the remaining propellant mass at time step $t$. Here, the term $\left(m_{\text{dry}} + m_{p,\text{rem}}\left(t\right)\right) \exp\left(-\frac{\Delta v_{\text{stk}} t_{r,\text{rep}}}{g_0 I_{\text{sp}}}\right)$ represents the sum of the satellite's dry mass and the remaining propellant mass after the time step $t_{r,\text{rep}}$. Additionally, the OOR service cost with refueling amount of {\color{black}$m_{p,\text{oor}}\left(k,t\right)$} is represented as {\color{black}follows:
\begin{equation}
    C_{\text{oor}}\left(k,t\right)=c_{v,\text{oor}} m_{p,\text{oor}}\left(k,t\right)+c_{f,\text{oor}},\quad k\in \mathcal{K}_t,
\end{equation}
which is derived from Eq \eqref{eqn: CER launch}. Since $m_{p,\text{oor}}$ is bounded by the service capacity ($M_{\text{oor}}$), the set of available values for $k$ ($\mathcal{K}_t$) is defined as follows:
\begin{equation}
    \mathcal{K}_t=\left\{k\in \mathbb{Z}_{++}:m_{p,\text{oor}}\left(k,t\right)\le M_{\text{oor}}\right\}.
\end{equation}}
Lastly, {\color{black}$u_{\text{oor}}\left(k,t;t'\right)$} is defined as{\color{black}
\begin{equation}\label{eqn: UtilityOfOOR}
    u_{\text{oor}}\left(k,t;t'\right)=u_{\text{oor},1}\left(t;t'\right)+u_{\text{oor},2}\left(k,t;t'\right)+u_{\text{oor},3}\left(k,t;t'\right),
\end{equation}}
{\color{black}where} $u_{\text{oor},1}$ represents the EA of the scenario where the OOR service fails; $u_{\text{oor},2}$ represents the EA of the scenario where the OOR service is successful but the replenished satellite experiences an in-orbit failure before reaching its operational lifetime; $u_{\text{oor},3}$ represents the EA of the scenario where the OOR service is successful and the replenished satellite operates until depletion without any in-orbit failure. The mathematical representation of each term is as follows:{\color{black}
\begin{equation}\label{eqn: UtilityOfOOR1-1}
    u_{\text{oor},1}\left(t;t'\right)=p_{\text{oor}} A/P\left(r_f,t_{\text{rep}}\right) EP_{1}\left(t;t'\right),
\end{equation}
\begin{align}
    u_{\text{oor,2}}\left(k,t;t'\right)=\left(1-p_{\text{oor}}\right)\sum_{t_1=t_{\text{rep}}+1}^{t_{\text{rep}}+k} P_{\text{iof}}\left(t_1+t;t'\right)A/P\left(r_f,t_1-1\right) &\Biggl(EP_{1}\left(t;t'\right)-C_{\text{oor}}\left(k,t\right)e^{-r_ft_{\text{rep}}\Delta T}\nonumber\\
    &+\sum_{t_2=t_{\text{rep}}+1}^{t_1}\mathbb{E}\left[R\left(t+t_2\right);t',t\right] e^{-r_ft_2\Delta T}\Biggr)\label{eqn: UtilityOfOOR1-2},
\end{align}
\begin{align}
    u_{\text{oor,3}}\left(k,t;t'\right)=\left(1-p_{\text{oor}}\right)\left(\sum_{\tau>t+t_{\text{rep}}+k}P_{\text{iof}}\left(\tau;t'\right)\right)A/P\left(r_f,t_{\text{rep}}+k\right)&\Biggl(EP_{1}\left(t;t'\right)-C_{\text{oor}}\left(k,t\right)e^{-r_ft_{\text{rep}}\Delta T}\nonumber\\
    &+\sum_{t_2=t_{\text{rep}}+1}^{t_{\text{rep}}+k}\mathbb{E}\left[R\left(t+t_2\right);t',t\right]e^{-r_ft_2\Delta T}\Biggr)\label{eqn: UtilityOfOOR1-3},
\end{align}
where $t'$ is the time step at which the currently operational satellite started its operation.}

\subsubsection{Decision 2}
In Decision 2, a range of alternatives is considered, each corresponding to a different amount of propellant mass required for refueling. Similar to the approach in Decision 1, the utility function for each alternative is defined based on the assumptions in Section \ref{sec: section 3}. In Decision 2, {\color{black}$u_{\text{oor}}$} denotes the utility of extending the operation by OOR by $k$ time steps. The formulation parallels Eq \eqref{eqn: UtilityOfOOR}, albeit with distinctions in the definitions of each term, as detailed below:{\color{black}
\begin{align}\label{eqn: UtilityOfOOR2-1}
    u_{\text{oor},1}\left(t;t'\right)=p_{\text{oor}} A/P\left(r_f,t_{\text{oor}}\right) EP_{2}\left(t;t'\right),
\end{align}
\begin{align}
    u_{\text{oor,2}}\left(k,t;t'\right)=\left(1-p_{\text{oor}}\right)\sum_{t_1=t_{\text{oor}}+1}^{t_{\text{oor}}+k-1} P_{\text{iof}}\left(t_1+t;t'\right) A/P\left(r_f,t_1-1\right) &\Biggl(EP_{2}\left(t;t'\right)-C_{\text{oor}}\left(k,t\right)e^{-r_ft_{\text{oor}}\Delta T}\nonumber\\
    &+\sum_{t_2=t_{\text{oor}}+1}^{t_1-1}\mathbb{E}\left[R\left(t+t_2\right);t',t\right] e^{-rt_2\Delta T}\Biggr),
    \label{eqn: UtilityOfOOR2-2}
\end{align}
\begin{align}
    u_{\text{oor,3}}\left(k,t;t'\right)=\left(1-p_{\text{oor}}\right)\left(\sum_{\tau>t+t_{\text{oor}}+k}P_{\text{iof}}\left(\tau;t'\right)\right)A/P\left(r_f,t_{\text{oor}}+k\right)&\Biggl(EP_2\left(t;t'\right)-C_{\text{oor}}\left(k,t\right)e^{-r_ft_{\text{oor}}}\Delta T\nonumber\\
    &+\sum_{t_2=t_{\text{oor}}+1}^{t_{\text{oor}}+k}\mathbb{E}\left[R\left(t+t_2\right);t',t\right]e^{-rt_2\Delta T}\Biggr).
    \label{eqn: UtilityOfOOR2-3}
\end{align}}

\section*{Acknowledgments}
This work was prepared at the Korea Advanced Institute of Science and Technology, Department of Aerospace Engineering, under a research grant from the National Research Foundation of Korea entitled ``Space Logistics Modeling and Demand Fulfillment Strategy Evaluation Framework.'' The authors thank the National Research Foundation of Korea for the support of this work.

\color{black}


\begin{thebibliography}{50}
\singlespacing
\bibitem{li2019orbit}
Li, W.-J., Cheng, D.-Y., Liu, X.-G., Wang, Y.-B., Shi, W.-H., Tang, Z.-X., Gao, F., Zeng, F.-M., Chai, H.-Y, Luo, W.-B., Cong, Q., and Gao, Z.-L., ``On-Orbit Service (OOS) of Spacecraft: A Review of Engineering Developments,'' \textit{Progress in Aerospace Sciences}, Vol. 108, Jul. 2019, pp. 32-120.
\\ \href{https://doi.org/10.1016/j.paerosci.2019.01.004}{https://doi.org/10.1016/j.paerosci.2019.01.004}

\bibitem{arney2021orbit}
Arney, D., Sutherland, R., Mulvaney, J., Steinkoenig, D., Stockdale, C., and Farley, M., ``On-Orbit Servicing, Assembly, and Manufacturing (OSAM) State of Play,'' NASA TR 20210022660, 2021.

\bibitem{davis2019orbit}
Davis, J. P., Mayberry, J. P., and Penn, J. P., ``Game Changer: On-Orbit Servicing,'' The Aerospace Corporation Center for Space Policy and Strategy, May 2019.
\\ \href{https://csps.aerospace.org/papers/game-changer-orbit-servicing}{https://csps.aerospace.org/papers/game-changer-orbit-servicing}

\bibitem{cavaciuti2022space}
Cavaciuti, A. J., Heying, J. H., and Davis, J., ``Game Changer: In-Space Servicing, Assembly, and Manufacturing for the New Space Economy,'' The Aerospace Corporation Center for Space Policy and Strategy, Jul. 2022.
\\ \href{https://csps.aerospace.org/papers/game-changer-space-servicing-assembly-and-manufacturing-new-space-economy}{https://csps.aerospace.org/papers/game-changer-space-servicing-assembly-and-manufacturing-new-space-economy}

\bibitem{orbitfab}
Orbit Fab, retrieved 10 June 2025.
\\ \href{https://www.orbitfab.com/}{https://www.orbitfab.com/}

\bibitem{hastings2016when}
Hastings, D. E., Putbrese, B. L., and La Tour, P. A., ``When Will On-Orbit Servicing be Part of the Space Enterprise?,'' \textit{Acta Astronautica}, Vol. 127, Oct.-Nov. 2016, pp. 655-666.
\\ \href{https://doi.org/10.1016/j.actaastro.2016.07.007}{https://doi.org/10.1016/j.actaastro.2016.07.007}

\bibitem{saleh2002spacecraft}
Saleh, J. H., Hastings, D. E., and Newman, D. J., ``Spacecraft Design Lifetime,'' \textit{Journal of Spacecraft and Rockets}, Vol. 39, No. 2, 2002, pp. 244-257.
\\ \href{https://doi.org/10.2514/2.3806}{https://doi.org/10.2514/2.3806}


\bibitem{saleh2004weaving}
Saleh, J. H., Hastings, D. E., and Newman, D. J., ``Weaving Time Into System Architecture: Satellite Cost per Operational Day and Optimal Design Lifetime,'' \textit{Acta Astronautica}, Vol. 54, No. 6, 2004, pp. 413-431.
\\ \href{https://doi.org/10.1016/S0094-5765(03)00161-9}{https://doi.org/10.1016/S0094-5765(03)00161-9}

\bibitem{saleh2006reduce}
Saleh, J. H., Torres-Padilla, J.-P., Hastings, D. E., and Newman, D. J., ``To Reduce or to Extend a Spacecraft Design Lifetime?,'' \textit{Journal of Spacecraft and Rockets}, Vol. 43, No. 1, 2006, pp. 207-217.
\\ \href{https://doi.org/10.2514/1.10991}{https://doi.org/10.2514/1.10991}

\bibitem{saleh2003flexibility}
Saleh, J. H., Lamassoure, E. S., Hastings, D. E., and Newman, D. J., ``Space Systems Flexibility Provided by On-Orbit Servicing: Part 1,'' \textit{Journal of Spacecraft and Rockets}, Vol. 40, No. 2, 2003, pp. 551-560.
\\ \href{https://doi.org/10.2514/2.3944}{https://doi.org/10.2514/2.3944}

\bibitem{saleh2002space}
Saleh, J. H., Lamassoure, E., and Hastings, D. E., ``Space Systems Flexibility Provided by On-Orbit Servicing: Part 1,'' \textit{Journal of Spacecraft and Rockets}, Vol. 39, No. 4, 2002, pp. 551-560.
\\ \href{https://doi.org/10.2514/2.3844}{https://doi.org/10.2514/2.3844}

\bibitem{lamassoure2002space}
Lamassoure, E., Saleh, J. H., and Hastings, D. E., ``Space Systems Flexibility Provided by On-Orbit Servicing: Part 2,'' \textit{Journal of Spacecraft and Rockets}, Vol. 39, No. 4, 2002, pp. 561-570.
\\ \href{https://doi.org/10.2514/2.3845}{https://doi.org/10.2514/2.3845}

\bibitem{joppin2006orbit}
Joppin, C., and Hastings, D. E., ``On-Orbit Upgrade and Repair: The Hubble Space Telescope Example,'' \textit{Journal of Spacecraft and Rockets}, Vol. 43, No. 3, 2006, pp. 614-625.
\\ \href{https://doi.org/10.2514/1.15496}{https://doi.org/10.2514/1.15496}

\bibitem{long2007orbit}
Long, A. M., Richards, M. G., and Hastings, D. E., ``On-Orbit Servicing: A New Value Proposition for Satellite Design and Operation,'' \textit{Journal of Spacecraft and Rockets}, Vol. 44, No. 4, 2007, pp. 964-976.
\\ \href{https://doi.org/10.2514/1.27117}{https://doi.org/10.2514/1.27117}

\bibitem{lamassoure2001framework}
Lamassoure, E. S., ``A Framework to Account for Flexibility in Modeling the Value of On-Orbit Servicing for Space Systems,'' Master's Thesis, Massachusetts Institute of Technology, Cambridge, MA, 2001.

\bibitem{saleh2002weaving}
Saleh, J. H., ``Weaving Time Into System Architecture: New Perspectives on Flexibility, Spacecraft Design Lifetime, and On-Orbit Servicing,'' Ph.D. Dissertation, Massachusetts Institute of Technology, Cambridge, MA, 2002.

\bibitem{joppin2004orbit}
Joppin, C., ``On-Orbit Servicing for Satellite Upgrades,'' Master's Thesis, Massachusetts Institute of Technology, Cambridge, MA, 2004.

\bibitem{long2005framework}
Long A. M., ``Framework for Evaluating Customer Value and the Feasibility of Servicing Architectures for On-Orbit Satellite Servicing,'' Master's Thesis, Massachusetts Institute of Technology, Cambridge, MA, 2005.

\bibitem{richards2006orbit}
Richards, M. G., ``On-Orbit Serviceability of Space System Architectures,'' Master's Thesis, Massachusetts Institute of Technology, MA, 2006.

\bibitem{mun2016real}
Mun, J., \textit{Real Options Analysis: Tools and Techniques for Valuing Strategic Investments and Decisions with Integrated Risk Management and Advanced Quantitative Decision Analytics}, $3^{\text{rd}}$ ed., ROV Press, Dublin, CA, 2016.


\bibitem{yao2013orbit}
Yao, W., Chen, X., Huang, Y., and van Tooren, M., ``On-Orbit Servicing System Assessment and Optimization Methods Based on Lifecycle Simulation Under Mixed Aleatory and Epistemic Uncertainties,'' \textit{Acta Astronautica}, Vol. 87, Jun.-July 2013, pp. 107-126.
\\ \href{https://doi.org/10.1016/j.actaastro.2013.02.005}{https://doi.org/10.1016/j.actaastro.2013.02.005}

\bibitem{liu2021economic}
Liu, Y., Zhao, Y., Tan, C., Liu, H., and Liu, Y., ``Economic Value Analysis of On-Orbit Servicing for Geosynchronous Communication Satellites,'' \textit{Acta Astronautica}, Vol. 180, Mar. 2021, pp. 176--188.
\\ \href{https://doi.org/10.1016/j.actaastro.2020.11.040}{https://doi.org/10.1016/j.actaastro.2020.11.040}

\bibitem{henry2020intelsat}
Henry, C., ``Intelsat-901 Satellite, With MEV-1 Servicer Attached, Resumes Service,'' Space News, 2020, retrieved 10 June 2025.
\\ \href{https://spacenews.com/intelsat-901-satellite-with-mev-1-servicer-attached-resumes-service/}{https://spacenews.com/intelsat-901-satellite-with-mev-1-servicer-attached-resumes-service/}

\bibitem{krebs}
Krebs, G. D., ``Intelsat 39,'' Gunter's Space Page, retrieved 10 June 2025.
\\ \href{https://space.skyrocket.de/doc_sdat/intelsat-39.htm}{https://space.skyrocket.de/doc\_sdat/intelsat-39.htm}

\bibitem{markowitz1952portfolio}
Markowitz, H., ``Portfolio Selection,'' \textit{The Journal of Finance}, Vol. 7, No. 1, 1952, pp. 77--91.
\\ \href{https://doi.org/10.2307/2975974}{https://doi.org/10.2307/2975974}

\bibitem{sharpe1966mutual}
Sharpe, W. F., ``Mutual Fund Performance,'' \textit{The Journal of Business}, Vol. 39, No. 1, 1966, pp. 119--138.
\\ \href{https://www.jstor.org/stable/2351741}{https://www.jstor.org/stable/2351741}

\bibitem{marler2004survey}
Marler, R. T., and Arora, J. S., ``Survey of Multi-Objective Optimization Methods for Engineering,'' \textit{Structural and Multidisciplinary Optimization}, Vol. 26, Mar. 2004, pp. 396-395.
\\ \href{https://doi.org/10.1007/s00158-003-0368-6}{https://doi.org/10.1007/s00158-003-0368-6}

\bibitem{newnan2017engineering}
Newman, D. G., Eschenbach, T. G., and Lavelle, J. P., \textit{Engineering Economic Analysis}, $13^{\text{th}}$ ed., Oxford University Press, Oxford, United Kingdom, 2017.

\bibitem{miles1980weighted}
Miles, J. A., and Ezzell, J. R., ``The Weighted Average Cost of Capital, Perfect Capital Markets, and Project Life: A Clarification,'' \textit{Journal of Financial and Quantitative Analysis}, Vol. 15, No. 3, 1980, pp. 719-730.
\\ \href{https://doi.org/10.2307/2330405}{https://doi.org/10.2307/2330405}

\bibitem{wertz2011space}
Wertz, J. R., Everett, D. F., and Puschell, J. J., \textit{Space Mission Engineering: The New SMAD}, Microcosm Press, Torrance, CA, 2011.

\bibitem{cpi}
``Consumer Price Index, 1913–'', Federal Reserve Bank of Minneapolis, 2025, retrieved 10 June 2025.
\\ \href{https://www.minneapolisfed.org/about-us/monetary-policy/inflation-calculator/consumer-price-index-1913-}{https://www.minneapolisfed.org/about-us/monetary-policy/inflation-calculator/consumer-price-index-1913-}

\bibitem{hull2022options}
Hull, J.  C., \textit{Options, Future, and Other Derivaties}, $11^{\text{th}}$, Pearson Education, London, United Kingdom, 2022.

\bibitem{saleh2011spacecraft}
Saleh, J. H., and Castet, J.-F., \textit{Spacecraft Reliability and Multi-State Failures: A Statistical Approach}, John Wiley \& Sons, Hoboken, NJ, 2011.
\\ \href{https://doi.org/10.1002/9781119994077}{https://doi.org/10.1002/9781119994077}

\bibitem{arianespace2021ariane}
Lagier, R., ``Ariane 6 User's Manual (Issue 2 Revision 0),'' Arianespace Corp., 2021, retrieved 10 June 2025.
\\ \href{https://www.arianespace.com/ariane-6/}{https://www.arianespace.com/ariane-6/}


\bibitem{rasmussen2006gaussian}
Rasmussen, C. E., and Williams, C. K. I., \textit{Gaussian Process for Machine Learning}, The MIT Press, Cambridge, MA, 2006.
\\ \href{https://doi.org/10.7551/mitpress/3206.001.0001}{https://doi.org/10.7551/mitpress/3206.001.0001}

\bibitem{montgomery2012design}
Montgomery, D. C., \textit{Design and Analysis of Experiments}, $8^{\text{th}}$ ed., John Wiley \& Sons, Hoboken, NJ, 2012.

\bibitem{rice2006mathematical}
Rice, J. A., \textit{Mathematical Statistics and Data Analysis}, $3^{\text{rd}}$ ed., Cengage, Independence, KY, 2006.

\bibitem{deb2002fast}
Deb, K., Pratap, A., Agarwal, S., and Meyarivan, T., ``A Fast and Elitist Multiobjective Genetic Algorithm: NSGA-II,'' \textit{IEEE Transactions on Evolutionary Computation}, Vol. 6, No. 2, 2009, pp. 182-197.
\\ \href{https://doi.org/10.1109/4235.996017}{https://doi.org/10.1109/4235.996017}

\bibitem{wall2024europe}
Wall, M., ``Europe's New Ariane 6 Rocket Launches on Long-Awaited Debut Mission (Video),'' Space.com, 2024, retrieved 10 June 2025.
\\ \href{https://www.space.com/europe-ariane-6-rocket-debut-launch}{https://www.space.com/europe-ariane-6-rocket-debut-launch}

\bibitem{satnews2024}
``USSF-51 for the United States Space Force's Space Systems Command,'' Satnews, 2024, retrieved 10 June 2025.
\\ \href{https://news.satnews.com/2024/07/02/ulas-atlas-v-rocket-to-launch-ussf-51-for-the-united-states-space-forces-space-systems-command/}{https://news.satnews.com/2024/07/02/ulas-atlas-v-rocket-to-launch-ussf-51-for-the-united-states-space-forces-space-systems-command/}

\bibitem{atlasv}
``Atlas V,'' United Launch Alliance, retrieved 10 June 2025.
\\ \href{https://www.ulalaunch.com/rockets/atlas-v}{https://www.ulalaunch.com/rockets/atlas-v}

\bibitem{falcon9}
``Falcon 9,'' SpaceX, retrieved 10 June 2025.
\\ \href{https://www.spacex.com/vehicles/falcon-9/}{https://www.spacex.com/vehicles/falcon-9/}

\bibitem{showalter2025progression}
Showalter, N. E., Hamel, J. M., de Weck, O. L., and Hasting, D. E., ``Progression of Satellite Refueling and Repositioning Technologies Through a Client-Servicer Perspective,'' \textit{Journal of Spacecraft and Rockets}, Published Online.
\\ \href{https://doi.org/10.2514/1.A36256}{https://doi.org/10.2514/1.A36256}

\bibitem{dujonchay2017quantification}
du Jonchay, T. S., and Ho, K., ``Quantification of the Responsiveness of On-Orbit Servicing Infrastructure for Modularized Earth-Orbiting Platforms,'' \textit{Acta Astronautica}, Vol. 132, 2017, pp. 192--203.
\\\href{https://doi.org/10.1016/j.actaastro.2016.12.021}{https://doi.org/10.1016/j.actaastro.2016.12.021}

\bibitem{ho2020semi}
Ho, K., Wang, H., DeTrempe, P. A., du Jonchay, T. S., and Tomita, K., ``Semi-Analytical Model for Design and Analysis of On-Orbit Servicing Architecture,'' \textit{Journal of Spacecraft and Rockets}, Vol. 57, No. 6, 2020, pp. 1129--1138.
\\\href{https://doi.org/10.2514/1.A34663}{https://doi.org/10.2514/1.A34663}

\bibitem{esa2004spt100}
``SPT-100,'' European Space Agency, retrieved 10 June 2025.
\\ \href{https://www.esa.int/esearch?q=SPT-100}{https://www.esa.int/esearch?q=SPT-100}

\bibitem{goebel2009evaluation}
Goebel, D. M., Polk, J. E., Sandler, I., Mikellides, I. G., and Brophy, J. R., ``Evaluation of 25-cm XIPS Thruster Life for Deep Space Mission Applications,'' University of Michigan, \textit{31st International Electric Propulsion Conference}, Paper IEPC-2009-152, Sept. 2009.

\end{thebibliography}
\end{document}